\def\no{\if01}
\def\iftwelvept{\no}

\def\ifusepdf{\no}
\def\ifpsfont{\no}

\iftwelvept
\documentclass[leqno,10pt]{amsart}
\else
\documentclass[leqno,12pt]{amsart}
\fi
\usepackage{amssymb}
\usepackage{amscd}
\usepackage{latexsym}
\usepackage{verbatim}
\usepackage{mathrsfs}
\usepackage[all]{xy}

\ifusepdf
\usepackage{hyperref}
\else\fi
\ifpsfont
\usepackage[T1]{fontenc}
\usepackage{times}
\else\fi

\iftwelvept
\else\fi

\theoremstyle{plain}
\newtheorem{Theorem}{Theorem}[section]

\newtheorem{Proposition}[Theorem]{Proposition}
\newtheorem{Lemma}[Theorem]{Lemma}
\newtheorem{Corollary}[Theorem]{Corollary}

\theoremstyle{definition}

\newtheorem{Definition}[Theorem]{Definition}
\newtheorem{Remark}[Theorem]{Remark}

\newcommand{\ZZ}{{\mathbb{Z}}}
\newcommand{\QQ}{{\mathbb{Q}}}
\newcommand{\RR}{{\mathbb{R}}}
\newcommand{\CC}{{\mathbb{C}}}

\newcommand{\NN}{{\mathbb{N}}}
\newcommand{\GG}{{\mathbb{G}}}

\newcommand{\XXX}{{\mathscr{X}}}
\newcommand{\YYY}{{\mathscr{Y}}}

\newcommand{\UUU}{{\mathscr{U}}}
\newcommand{\VVV}{{\mathscr{V}}}
\newcommand{\DDD}{{\mathscr{D}}}

\newcommand{\OOO}{{\mathscr{O}}}
\newcommand{\WWW}{{\mathscr{W}}}
\newcommand{\AAA}{{\mathscr{A}}}
\newcommand{\SIGMA}{{\mbox{\boldmath$\Sigma$}}}
\newcommand{\OO}{{\mathcal{O}}}

\newcommand{\LLL}{{\mathcal{L}}}
\newcommand{\MMM}{\mathcal{M}}
\newcommand{\NNN}{\mathcal{N}}
\newcommand{\SSS}{\mathcal{S}}

\newcommand{\Hom}{\textup{Hom}\,}
\newcommand{\extup}{\scriptsize{\textup{can}}}
\newcommand{\xtup}{\scriptsize{\textup{red}}}

\newcommand{\Spec}{\textup{Spec}\,}

\newcommand{\rank}{\textup{rk}\,}

\newcommand{\Aut}{\textup{Aut}\,}

\newcommand{\XX}{\mathcal{X}}
\newcommand{\YY}{\mathcal{Y}}

\newcommand{\Proof}{{\sl Proof.}\quad}
\newcommand{\QED}{{\unskip\nobreak\hfil\penalty50\quad\null\nobreak\hfil
{$\Box$}\parfillskip0pt\finalhyphendemerits0\par\medskip}}

\begin{document}
\title{The category of toric stacks}
\author{Isamu Iwanari}
\renewcommand{\thefootnote}{\fnsymbol{footnote}}
\footnote[0]{2000\textit{ Mathematics Subject Classification}.
 14M25, 14A20.}

\address{Research Institute for the Mathematical Sciences,
Kyoto University, Kyoto, 606-8502, Japan}
\email{iwanari@kurims.kyoto-u.ac.jp }
\begin{abstract}
In this paper, we prove that there exists an equivalence between
2-category of smooth Deligne-Mumford stacks with torus-embeddings and actions,
and the 1-category of stacky fans.
For this purpose, we obtain two main results.
The first is to investigate a combinatorial
aspect of the 2-category of toric algebraic stacks defined in \cite{I2}.
We establish an equivalence between the 2-category of toric algebraic stacks and
the 1-category of stacky fans.
The second is to give a geometric characterization theorem
for toric algebraic stacks.
\end{abstract}

\maketitle

\section{Introduction and main results}

The equivalence between the category of toric varieties
and that of fans is a fundamental theorem of toric varieties,
and it gives a fruitful bridge between algebraic geometry and combinatorics.
It is also useful in various stages; for example, 
the most typical and beautiful use is toric
minimal model program (cf. \cite{R}).
Since 
simplicial toric varieties have quotient singularities
in characteristic zero, it is a natural problem to find such an equivalence in the stack-theoretic
context. Let $k$ be an algebraically closed base field of characteristic zero.
Consider a triple
$(\mathcal{X}, \iota:\GG_m^d\hookrightarrow \mathcal{X}, a:\mathcal{X}\times\GG_m^d\to \mathcal{X})$, where 
$\mathcal{X}$ is a smooth Deligne-Mumford stack
of finite type and  separated over $k$, that satisfies the followings:
\begin{enumerate}
\renewcommand{\labelenumi}{(\roman{enumi})}

\item The morphism $\iota:\GG_m^d\hookrightarrow \mathcal{X}$ is an open immersion
identifying $\GG_m^d$ with a dense open substack of $\XX$. (We shall refer to $\GG_m^d\hookrightarrow \mathcal{X}$
to a torus-embedding.)

\item The morphism $a:\mathcal{X}\times\GG_m^d\to \mathcal{X}$
is an action of $\GG_m^d$ on $\mathcal{X}$,
which is an extension of the action
$\GG_m^d$ on itself. (We shall refer to it as a torus action.)

\item The coarse moduli space $X$ for $\mathcal{X}$ is a
scheme.

\end{enumerate}

We shall refer to such a triple as a {\it toric triple}.
Note that if $\XX$ is a scheme, then $\XX$ is a smooth
toric variety.
A 1-morphism of toric triples
\[
(\mathcal{X}, \iota:\GG_m^d\hookrightarrow \mathcal{X}, a:\mathcal{X}\times\GG_m^d\to \mathcal{X}) \to (\mathcal{X}', \iota':\GG_m^{d'}\hookrightarrow \mathcal{X}', a':\mathcal{X}'\times\GG_m^{d'}\to \mathcal{X}')
\]
is a morphism $f:\mathcal{X}\to \mathcal{X}'$
such that
the restriction of $f$ to $\GG_m^d$
induces a morphism $\GG_m^d\to \GG_m^{d'}$
of group $k$-schemes and the diagram
\[
\xymatrix@R=6mm @C=17mm{
\mathcal{X}\times\GG_m^d \ar[r]^{f\times (f|_{\GG_m^d})} \ar[d]^a & \mathcal{X}'\times\GG_m^{d'} \ar[d]^{a'} \\
\mathcal{X}\ar[r]^f &\mathcal{X}'  \\
}
\]
commutes in the 2-categorical sense.
A 2-isomorphism $g:f_1\to f_2$ is an isomorphism of 1-morphisms.

Our main goal is the following:

\begin{Theorem}
\label{main3}
There exists an equivalence between the 2-category of toric triples
and the 1-category of stacky fans. $($See Definition~\ref{stfan}
for the definition of stacky fans.$)$

\end{Theorem}

The non-singular fans form a full-subcategory of the category of stacky fans.
Thus our equivalence includes
the classical equivalence between
smooth toric varieties and non-singular fans (See Remark~\ref{pict}).

To obtain Theorem~\ref{main3}, we need to consider the following two problems:

\begin{enumerate}
\renewcommand{\labelenumi}{(\roman{enumi})}

\item Construction of toric triples $\XXX_{(\Sigma,\Sigma^0)}$
associated to a stacky fan $(\Sigma,\Sigma^0)$
(we shall call it the associated toric algebraic stack),
and the establishment of an equivalence between the groupoid category $\Hom(\XXX_{(\Sigma_1,\Sigma_1^0)},\XXX_{(\Sigma_2,\Sigma_2^0)})$ and the discrete category associated to
the set of morphisms $\Hom((\Sigma_1,\Sigma_1^0),(\Sigma_2,\Sigma_2^0))$.

\item Geometric characterization of toric algebraic stacks
associated to stacky fans.

\end{enumerate}

\vspace{2mm}

For the first problem, the construction was given
in \cite{I2} (see also \cite{I}).
In loc. cit., given a stacky fan $(\Sigma,\Sigma^0)$ we defined the associated toric 
algebraic stack $\XXX_{(\Sigma,\Sigma^0)}$ by means of logarithmic geometry.
In characteristic zero it is a smooth
Deligne-Mumford stack and has a natural torus embedding and
a torus action, i.e., a toric triple.
Let us denote by $\mathfrak{Torst}$ 
the 2-category whose objects
are toric algebraic stacks associated to stacky fans (cf. Section 2). A 1-morphism of two toric 
algebraic stacks 
$f:\XXX_{(\Sigma,\Sigma^0)}\to \XXX_{(\Delta,\Delta^0)}$
in $\mathfrak{Torst}$
is
a torus-equivariant $1$-morphism (cf. Definition~\ref{stfan},~\ref{torst},~\ref{torusequiv}).
A 2-morphism $g:f_1\to f_2$ is an isomorphism of 1-morphisms.
Then the following is our answer to the first problem.

\begin{Theorem}
\label{main}
We work over a field of characteristic zero.
There exists an equivalence of 2-categories
\[
 \Phi:  \mathfrak{Torst}\stackrel{\sim}{\longrightarrow} (\textup{1-Category\ of\ Stacky\ fans})
\]
which makes the following diagram
\[
\begin{CD}
\mathfrak{Torst} @>{\Phi}>>  (\textup{Category\ of\ Stacky\ fans})\\
@V{c}VV            @VVV   \\
{\textup{Simtoric}} @>{\sim}>> (\textup{Category\ of\ simplicial\ fans})
\end{CD}
\]
commutative. Here ${\textup{Simtoric}}$ is the category of simplicial
toric varieties $($morphisms in ${\textup{Simtor}}$
are torus-equivariant ones$)$, and $c$ is the natural functor
which sends the toric algebraic stack $\XXX_{(\Sigma,\Sigma^0)}$
associated to a stacky fan $(\Sigma,\Sigma^0)$
to the
toric variety $X_{\Sigma}$ $($cf. Remark~\ref{degefun}
and Remark~\ref{des}$)$.
The 1-category of stacky fans is regarded as a 2-category.
\end{Theorem}

In Theorem~\ref{main} the difficult point is to show that the groupoid of torus-equivariant $1$-morphisms of toric algebraic stacks is equivalent to the discrete category
of the set of the morphisms of stacky fans.
For two algebraic stacks $\XX$ and $\YY$, the classification of 1-morphisms
$\XX\to \YY$ and their 2-isomorphisms is a hard problem even if $\XX$ and 
$\YY$ have very explicit groupoid presentations.
Groupoid presentations of stacks are ill-suited to considering such problems.
Our idea to overcome this difficulty
is to use the modular interpretation of toric algebraic stacks
in terms of {\it logarithmic geometry} (cf. Section 2).
The notion of a certain type of resolutions of monoids (called $(\Sigma,\Sigma^0)$-free
resolutions) plays a role similar to monoid-algebras arising from cones
in classical toric geometry,
and by virtue of it we obtain Theorem~\ref{main} by reducing it to a certain
problem on log structures on schemes.
As a corollary of Theorem~\ref{main}, we show also that every toric algebraic stack
admits a smooth torus-equivariant cover by a smooth toric variety
(cf. Corollary~\ref{cor}).

\vspace{1mm}

The second problem is the geometric characterization.
Remembering that we have the {\it geometric}
characterization of toric varieties (cf. \cite{KKMS}) we wish to
have a similar characterization of toric algebraic stacks.
Our geometric characterization of toric algebraic stacks is:

\begin{Theorem}
\label{main2}
Assume that the base field $k$ is algebraically closed in characteristic zero.
Let $(\mathcal{X}, \iota:\GG_m^d\hookrightarrow \mathcal{X}, a:\mathcal{X}\times\GG_m^d\to \mathcal{X})$ be a toric triple over $k$.
Then there exist a stacky fan $(\Sigma,\Sigma^0)$
and an isomorphism of stacks
\[
\Phi:\mathcal{X}\stackrel{\sim}{\longrightarrow} \XXX_{(\Sigma,\Sigma^0)}
\]
over $k$,
and it satisfies the following additional properties:
\begin{enumerate}

\item The restriction of $\Phi$ to $\GG_m^d\subset \XX$
induces an isomorphism $\Phi_0:\GG_m^d\stackrel{\sim}{\to} \Spec k[M]\subset \XXX_{(\Sigma,\Sigma^0)}$
of group $k$-schemes.
Here $N=\ZZ^d$,
$M=\Hom_{\ZZ}(N,\ZZ)$, $\Sigma$ is a fan in $N\otimes_{\ZZ}\RR$,
and $\Spec k[M]\hookrightarrow \XXX_{(\Sigma,\Sigma^0)}$
is the natural torus embedding $($cf. Definition~\ref{torst}$)$.

\item The diagram
\[
\xymatrix{
\XX\times \GG_m^d \ar[d]^{m}\ar[r]^(0.4){\Phi\times\Phi_0} & \XXX_{(\Sigma,\Sigma^0)}\times\Spec k[M] \ar[d]^{a_{(\Sigma,\Sigma^0)}}\\
\XX\ar[r]^{\Phi} & \XXX_{(\Sigma,\Sigma^0)} \\
}
\]
commutes in the 2-categorical sense, where $a_{(\Sigma,\Sigma^0)}:\XXX_{(\Sigma,\Sigma^0)}\times\Spec k[M]
\to \XXX_{(\Sigma,\Sigma^0)}$ is the torus action functor
(cf. Section 2.2).
\end{enumerate}
Moreover, such a stacky fan is unique up to isomorphisms.
\end{Theorem}

The essential point in the proof of Theorem~\ref{main2} is
a study of (\'etale) local structures of the coarse moduli map
$\XX\to X$. We first show that $X$ is a toric variety
and then determine the local structure of $\XX\to X$
by applying logarithmic Nagata-Zariski purity Theorem due to
S. Mochizuki and K. Kato (independently proven).

\vspace{1mm}

It is natural and interesting to consider a generalization of our work
to positive characteristics.
Unfortunately, our proof does not seem applicable to the case of positive
characteristics. For example, we use the assumption of characteristic zero
for the application of log purity theorem.
Furthermore, in positive characteristics, toric algebraic stacks
defined in \cite{I2} are not necessarily Deligne-Mumford stacks.
It happens to be an Artin stack. Thus in such a generalization
the formulation should be modified. (See Remark~\ref{semifinalrem})

\vspace{1mm}

Informally, Theorem~\ref{main3} implies that the geometry of toric triples could be
encoded by the combinatorics of stacky fans.
It would be interesting to investigate the geometric invariants
of toric triples in the view of stacky fans.
In this direction, for instance, in the subsequent paper \cite{I3}
we have shown the relationship between integral Chow rings
of toric triples and (classical) Stanley-Reisner rings, in which
a non-scheme-theoretic phenomenon arises.
In another direction, it seems that Theorem~\ref{main}
has a nice place in the study of toric minimal
model program from a stack-theoretic and derived categorical
viewpoint.

\smallskip

The paper is organized as follows.
In Section 2 we recall basic definitions concerning toric algebraic stacks
and some results which we use in Section 3 and 4.
In Section 3 we present the proof of Theorem~\ref{main} and its corollaries.
In Section 4 we give the proof of Theorem~\ref{main2}.
Finally, in Section 5 we discuss the relationships
with the work \cite{BCS} of Borisov-Chen-Smith modeling the quotient
construction by Cox \cite{C}, and the recent works \cite{FMN}, \cite{P}
by Fantechi-Mann-Nironi and Perroni.

\smallskip

We systematically use the language of logarithmic geometry, and assume that readers are familiar with it at the level of the paper \cite{log}.

{\bf Notations And Conventions}

\smallskip
\begin{enumerate}
\renewcommand{\labelenumi}{(\theenumi)}

 \item We fix a Grothendieck universe $\mathscr{U}$
 with $\{ 0,1,2,3,...\}\in \UUU$ where 
 $\{ 0,1,2,3,...\}$ is the set of all finite ordinals.
 We consider only monoids, groups, rings, schemes
 and log schemes which belong to $\mathscr{U}$.

\item A {\it variety} is
a geometrically integral scheme of finite type and separated over a field.

\item {\it toric varieties}\,:
Let $N\cong \ZZ^d$ be a lattice of rank $d$
and $M=\Hom_{\ZZ}(N,\ZZ)$.
If $\Sigma$ is a fan in $N_{\RR}=N\otimes_{\ZZ}\RR$,
we denote the associated toric variety by $X_{\Sigma}$.
We usually write $i_{\Sigma}:T_{\Sigma}:=\Spec k[M]\hookrightarrow X_{\Sigma}$
for the torus embedding.
If $k$ is an algebraically closed field, then
by applying Sumihiro's Theorem \cite[Corollary 3.11]{S} just as in
\cite[Chap. 1]{KKMS}
we have the following geometric characterization of toric varieties:
Let $X$ be a normal variety which contains an algebraic torus ($= \GG_m^d$)
as a dense open subset. Suppose that the action of $\GG_m^d$ on itself
extends to an action of $\GG_m^d$ on $X$.
Then there exist a fan $\Sigma$ and
an equivariant isomorphism $X\cong X_{\Sigma}$.

\item {\it logarithmic geometry}\,:
All monoids will be assumed to be commutative with unit.
For a monoid $P$, we denote by $P^{gp}$ the Grothendieck group of $P$.
A monoid $P$ is said to be {\it sharp} if whenever $p+p'=0$ for $p,p'\in P$,
then $p=p'=0$.
For a fine sharp monoid $P$, an element $p\in P$ is said to
be {\it irreducible} if whenever
$p=q+r$ for $q,r\in P$, then either $q=0$ or $r=0$.
In this paper, a log structure on a scheme
$X$ means a log structure (in the sense of Fontaine-Illusie (\cite{log}))
on the \'etale site $X_{et}$.
We usually denote simply by $\MMM$ a log structure $\alpha:\MMM\to \OOO_{X}$
on $X$,
and denote by $\bar{\MMM}$ the sheaf $\MMM/\OOO_{X}^*$.
Let $R$ be a ring.
For a fine monoid $P$, the canonical log structure, denoted by $\MMM_{P}$,
on $\Spec R[P]$ is the log structure associated to the natural injective map
$P\to R[P]$.
If there is a homomorphism of monoids $P\to R$ (here we regard $R$
as a monoid under multiplication),
we denote by $\Spec (P\to R)$
the log scheme with underlying scheme $\Spec R$
and the log structure associated to $P\to R$.
For a toric variety $X_{\Sigma}$, we denote by $\MMM_{\Sigma}$
the fine log structure $\OOO_{X_{\Sigma}}\cap i_{\Sigma*}\OOO_{T_{\Sigma}}^*
\hookrightarrow \OOO_{X_{\Sigma}}$ on $X_{\Sigma}$.
We shall refer to this log structure as the {\it canonical log structure}
on $X_{\Sigma}$.
We refer to \cite{I2} for the further generalities and notations
concerning toric varieties, monoids and log schemes, which are required
in what follows.

 \item {\it algebraic stacks}\,: We follow the conventions of (\cite{LM}).
 For a diagram $X\stackrel{a}{\to} Y\stackrel{b}{\leftarrow} Z$,
 we denote by $X\times_{a,Y,b}Z$ the fiber product (and we often omit
 $a$ and $b$ if no confusion seems to likely arise).
 Let us review some facts on coarse moduli spaces of algebraic stacks.
 Let $\XX$ be an algebraic stack over a scheme $S$.
 A {\it coarse moduli space (or map)}
 for $\XX$ is a morphism $\pi:\XX\to X$ to an algebraic space
 over $S$
 such that (i) $\pi$ is universal among morphisms from $\XX$ to
 algebraic spaces over $S$, and (ii) for every algebraically closed
 $S$-field $K$ the map $[\XX(K)]\to X(K)$
 is bijective where $[\XX(K)]$ denotes the set of
 isomorphism classes of objects in the small category $\XX(K)$.
 The fundamental existence theorem on coarse
 moduli spaces (to which we refer as Keel-Mori Theorem \cite{KM})
 is stated as follows
 (the following version is enough for our purpose):
 Let $k$ be a field of characteristic zero.
 Let $\XX$ be an algebraic stack of finite type over $k$ with
 finite diagonal.
 Then there exists a coarse moduli space $\pi:\XX\to X$
 where $X$ is of finite type and separated over $k$, and
 it satisfies the additional properties:
 (a) $\pi$ is proper, quasi-finite and surjective,
 (b) For any morphism $X'\to X$ of algebraic spaces over $k$,
 $\XX\times_{X}X'\to X'$ is a coarse moduli space
 (cf. \cite[Lemma 2.3.3, Lemma 2.2.2]{AV}).

\end{enumerate}

{\it Acknowledgements} I would like to thanks Yuichiro Hoshi
for explaining to me basic facts on Kummer log \'etale covers and
log fundamental groups, and Prof. Fumiharu Kato
for his valuable comments.
I also want to thank the referee for his/her helpful comments.
I would like to thank
Institut de Math\'ematiques de Jussieu
for the hospitality during the stay where
a part of this work was done.

\section{Preliminaries}

In this Section, we will recall the basic definitions and properties
(\cite{I2}) concerning toric algebraic stacks and stacky fans.
We fix a base field $k$ of characteristic zero.

\subsection{Definitions}
In this paper, all fans are assumed to be {\it finite},
though the theory in \cite{I2} works also in the case of {\it infinite fans}.
For a fan $\Sigma$, we denote by $\Sigma(1)$ the set of rays.

\begin{Definition}
\label{stfan}
Let $N\cong \ZZ^d$ be a lattice of rank $d$
and $M=\Hom_{\ZZ}(N,\ZZ)$ the dual lattice.
A stacky fan is a pair $(\Sigma, \Sigma^0)$,
where $\Sigma$ is a simplicial fan
in $N_{\RR}=N\otimes_{\ZZ}\RR$, and
$\Sigma^0$ is a subset of $|\Sigma|\cap N$,
called the {\it free-net of $\Sigma$}, which has
the following property $(\spadesuit)$:

$(\spadesuit)$ For any cone $\sigma$ in $\Sigma$, $\sigma\cap \Sigma^0$ is
a submonoid of $\sigma\cap N$ which is
isomorphic to $\NN^{\dim \sigma}$,
such that for any element $e \in \sigma\cap N$
there exists a positive integer $n$
such that $n\cdot e \in \sigma\cap \Sigma^0$.

A morphism $f:(\Sigma\ \textup{in}\ N\otimes_{\ZZ}\RR, \Sigma^0)\to (\Delta\ \textup{in}\ N'\otimes_{\ZZ}\RR,\Delta^0)$ 
is a homomorphism of $\ZZ$-modules $f:N \to N'$ which satisfies
the following properties:

$\bullet$ For any cone $\sigma$ in $\Sigma$, there exists a cone $\tau$ in $\Delta$ such that $f\otimes_{\ZZ}\RR(\sigma)\subset \tau$.

$\bullet$ $f(\Sigma^0) \subset \Delta^0$.
\end{Definition}
There exists a natural forgetting functor
\[
(\textup{Category\ of\ stacky\ fans})\to
(\textup{Category\ of\ simplicial\ fans}),\ (\Sigma,\Sigma^0)\mapsto \Sigma.
\]
It is essentially surjective, but {\it not } fully faithful.
Given a stacky fan $(\Sigma,\Sigma^0)$ and a ray $\rho$ in $\Sigma(1)$,
the initial point $P_{\rho}$ of $\rho\cap \Sigma^0$ is said to be
the {\it generator} of $\Sigma^0$ on $\rho$.
Let $Q_{\rho}$ be the first point of $\rho\cap N$ and let $n_{\rho}$
be the natural number such that $n_{\rho}\cdot Q_{\rho}=P_{\rho}$.
Then the number $n_{\rho}$ is said to be the {\it level} of $\Sigma^0$
on $\rho$. Notice that $\Sigma^0$ is completely determined by the levels
of $\Sigma^0$ on rays of $\Sigma$.
Each simplicial fan $\Sigma$ has the canonical free-net
$\Sigma^0_{\extup}$ whose level on every ray in $\Sigma$
is one.

If $\Sigma$ and $\Delta$ are non-singular fans, then
a usual
morphism of fans $\Sigma\to \Delta$ amounts to a morphism of stacky fans $(\Sigma,\Sigma_{\extup}^0)\to (\Delta,\Delta_{\extup})$.
Namely, the category of non-singular fans is a full-subcategory of
the category of stacky fans.

Let us give an example. Let $\sigma$ be a 2-dimensional cone in $(\ZZ\cdot e_1\oplus \ZZ\cdot e_2)\otimes_{\ZZ}\RR=\RR^{\oplus 2}$, that is generated by
$e_1$ and $e_1+2e_2$.
Let $\sigma^0$ be a free submonoid of $\sigma\cap(\ZZ\cdot e_1\oplus \ZZ\cdot e_2)$ that is generated by $2e_1$ and $e_1+2e_2$.
Note $\sigma^0\cong \NN^{\oplus2}$.
Then $(\sigma, \sigma^0)$ forms a stacky fan.
The level of $\sigma^0$ on the ray $\RR_{\ge0}\cdot e_1$ (resp. $\RR_{\ge0}\cdot (e_1+2e_2)$) is 2 (resp. 1).

\vspace{2mm}

Let $P$ be a monoid and let $S\subset P$ be a submonoid.
We say that $S$ is {\it close to} $P$
if for any element
$e$ in $P$ there exists a positive integer $n$ such that
$n\cdot e$ lies in $S$.
The monoid $P$ is said to be {\it toric} if $P$ is a fine, saturated
 and torsion-free monoid.

Let $P$ be a toric sharp monoid and
$d$ the rank of $P^{gp}$. 
A toric sharp monoid $P$ is said to be
{\it simplicially toric} if there
exists a submonoid $Q$ of $P$ generated by $d$ elements
such that $Q$ is close to $P$.
\begin{Definition}
\label{mini}
Let $P$ be a simplicially toric sharp monoid and $d$ the rank
of $P^{gp}$.
The {\it minimal free resolution} of $P$
is an injective homomorphism of monoids
\[
i:P \longrightarrow F
\]
with $F\cong\NN^d$, which has the
following properties:
\begin{enumerate}
\renewcommand{\labelenumi}{(\theenumi)}

\item  The submonoid $i(P)$ is close to $F$.

\item For any injective homomorphism $j:P \to G$
such that $j(P)$ is  close to $G$ and $G\cong \NN^d$,
there exists a unique homomorphism $\phi:F\to G$
such that $j=\phi\circ i$.
\end{enumerate}
\end{Definition}

We remark that by \cite[Proposition 2.4]{I2} or Lemma~\ref{constfreeres}
there exists a unique minimal free resolution for any simplicially toric
sharp
monoid. Next we recall the definition of toric algebraic stacks 
(\cite{I2}). Just after Remark~\ref{degefun}, we give another
definition of toric algebraic stacks, which is more direct
presentation in the terms of logarithmic geometry.

\begin{Definition}
\label{torst}
The {\it toric algebraic stack} associated to a stacky fan $(\Sigma\ \textup{in}\ N\otimes_{\ZZ}\RR,\Sigma^0)$
is a stack $\XXX_{(\Sigma,\Sigma^0)}$ over the category of $k$-schemes
whose objects over a $k$-scheme $X$ are triples
$(\pi:\SSS \to \OO_X,\alpha:\MMM\to \OO_X,\eta:\SSS\to \MMM)$
such that:

\begin{enumerate}
\renewcommand{\labelenumi}{(\theenumi)}

\item $\SSS$ is an \'etale sheaf of submonoids of the constant
sheaf $M$ on $X$ determined by $M=\Hom_{\ZZ} (N,\ZZ)$ such
that for every point $x \in X$, $\SSS_{x}\cong \SSS_{\bar{x}}$.
Here $\SSS_{x}$ (resp. $\SSS_{\bar{x}}$) denotes the Zariski
(resp. \'etale) stalk.

\item $\pi:\SSS\to \OO_X$ is a map of monoids where $\OO_X$ is a monoid
under multiplication.

\item For $s \in \SSS, \pi(s)$ is invertible if and only if
$s$ is invertible.

\item For each point $x \in X$, there exists some $\sigma\in \Sigma$
such that $\SSS_{\bar{x}}= \sigma^{\vee}\cap M$.

\item $\alpha:\MMM\to\OO_X$ is a fine log structure on $X$.

\item $\eta:\SSS \to \MMM$ is a homomorphism of sheaves of monoids
such that $\pi=\alpha\circ \eta$,
and for each geometric point $\bar{x}$ on $X$, $\bar{\eta}:\bar{\SSS}_{\bar{x}}
=(\SSS/(\textup{invertible\ elements}))_{\bar{x}}\to
\bar{\MMM}_{\bar{x}}$ is isomorphic to the composite
\[
\bar{\SSS}_{\bar{x}} \stackrel{r}{\hookrightarrow} F\stackrel{t}{\hookrightarrow} F,
\]
where $r$ is the minimal free resolution of $\bar{\SSS}_{\bar{x}}$
and $t$ is defined as follows:
\end{enumerate}
Each irreducible element of $F$ canonically corresponds to
a ray in $\Sigma$ (See Lemma~\ref{irr} below).
Let us denote by $e_{\rho}$ the irreducible element of $F$
which corresponds to a ray $\rho$.
Then define $t:F\to F$
by $e_{\rho} \mapsto n_{\rho}\cdot e_{\rho}$
where $n_{\rho}$ is the level of $\Sigma^0$ on $\rho$.
We shall refer to $t\circ r:\bar{\SSS}_{\bar{x}}\to F$
as the {\it$(\Sigma,\Sigma^0)$-free resolution at $\bar{x}$}
(or {\it $(\Sigma,\Sigma^0)$-free resolution of $\bar{\SSS}_{\bar{x}}=\SSS_{\bar{x}}/(\textup{invertible elements}$})).

A set of morphisms from
$(\pi:\SSS \to \OO_X,\alpha:\MMM\to \OO_X,\eta:\SSS\to \MMM)$
to $(\pi':\SSS' \to \OO_X,\alpha':\MMM'\to \OO_X,\eta':\SSS'\to \MMM')$
over $X$
is the set of
isomorphisms of log structures $\phi:\MMM\to\MMM'$
such that $\phi\circ \eta=\eta':\SSS=\SSS'\to \MMM'$
if $(\SSS,\pi)=(\SSS',\pi')$
and
is an empty set if $(\SSS,\pi)\ne(\SSS',\pi')$.
With the natural notion of pullbacks, $\XXX_{(\Sigma,\Sigma^0)}$ is 
a fibered category.

By \cite[Theorem on page 10]{AMRT},
$\Hom_{k\textup{-schemes}}(X,X_{\Sigma})\cong \{$ all pair $(\SSS,\pi)$ on $X$ satisfying
(1),(2),(3),(4) $\}$.
Therefore there exists a natural functor $\pi_{(\Sigma,\Sigma^0)}:\XXX_{(\Sigma,\Sigma^0)} \longrightarrow X_{\Sigma}$
which simply forgets the data $\alpha:\MMM\to \OO_X$ and $\eta:\SSS\to \MMM$.
Moreover $\alpha:\MMM\to \OO_X$ and $\eta:\SSS\to \MMM$ are morphisms of the \'etale sheaves and thus $\XXX_{(\Sigma,\Sigma^0)}$ is a stack with respect to the
\'etale topology.
Objects of the form
$(\pi:M\to \OOO_X,\OOO_X^*\hookrightarrow \OOO_X,\pi:M\to \OOO_X^*)$
determine a full sub-category of $\XXX_{(\Sigma,\Sigma^0)}$,
i.e., the natural inclusion $i_{(\Sigma,\Sigma^0)}:T_{\Sigma}=
\Spec k[M]\hookrightarrow \XXX_{(\Sigma,\Sigma^0)}$.
This commutes with the torus-embedding
$i_{\Sigma}:T_{\Sigma}\hookrightarrow X_{\Sigma}$.
\end{Definition}

\begin{Lemma}
\label{irr}
With notation in Definition~\ref{torst},
let $e$ be an irreducible element
in $F$ and let $n$ be a positive integer such that
$n\cdot e\in r(\bar{\SSS}_{\bar{x}})$.
Let $m\in\SSS_{\bar{x}}$ be a lifting of $n\cdot e$.
Suppose that $\SSS_{\bar{x}}= \sigma^{\vee}\cap M\subset M$.
Then there exists a unique ray $\rho\in \sigma(1)$ such that
$\langle m,\zeta_{\rho}\rangle>0$,
where $\zeta_{\rho}$
is the first lattice point of $\rho$,
and $\langle \bullet,\bullet\rangle$ is the dual pairing.
It does not depend on the choice of liftings.
Moreover this correspondence defines a natural injective map
\[
\{ \textup{Irreducible\ elements\ of } F \} \to \Sigma(1).
\]
\end{Lemma}
\Proof
This is fairly elementary (and follows from \cite{I2}), but we will give the proof for the completeness.
Since the kernel of $\SSS_{\bar{x}}\to \bar{\SSS}_{\bar{x}}$
is $\sigma^{\perp}\cap M$, thus $\langle m,v_{\rho}\rangle$ does not depend upon the choice of liftings $m$.
Taking a splitting $N\cong N'\oplus N''$ such that
$\sigma\cong \sigma'\oplus \{ 0\}\subset N'_{\RR}\oplus N''_{\RR}$ where $\sigma'$ is a full-dimensional cone in
$N'_{\RR}$, we may and will assume that
$\sigma$ is a full-dimensional cone,
i.e., $\sigma^{\vee}\cap M$ is sharp.
Let $\iota:\sigma^{\vee}\cap M\hookrightarrow \sigma^{\vee}$ be
the natural inclusion
and $r:\sigma^{\vee}\cap M\hookrightarrow F$
the minimal free resolution.
Then there exists a unique injective homomorphism
$i:F\to \sigma^{\vee}$ such that $i\circ r=\iota$.
By this embedding, we regard $F$ as a submonoid of $\sigma^{\vee}$.
Since $r:\sigma^{\vee}\cap M\hookrightarrow F\subset \sigma^{\vee}$
is the minimal free resolution and $\sigma^{\vee}$ is
a simplicial cone,
thus for each ray $\rho\in \sigma^{\vee}(1)$,
the initial point of $\rho\cap F$ is an irreducible
element of $F$. Since $\rank F^{gp}=\rank (\sigma^{\vee}\cap M)^{gp}=\dim \sigma^{\vee}=\dim \sigma$, thus each irreducible element of $F$ lies on
one of rays of $\sigma^{\vee}$.
It gives rise to a natural bijective map from
the set of irreducible elements of $F$ to $\sigma^{\vee}(1)$. 
Since $\sigma$ and $\sigma^{\vee}$ are simplicial, we have a natural bijective map $\sigma^{\vee}(1)\to \sigma(1);
\ \rho\mapsto \rho^{\star}$,
where $\rho^{\star}$ is the unique ray which
does not lie in $\rho^{\perp}$.
Therefore the composite map from
the set of irreducible elements of $F$ to $\sigma(1)$
is a bijective map. Hence it follows our claim.
\QED

\begin{Remark}
\begin{enumerate}
\item The above definition works over arbitrary base schemes.

\item If $\Sigma$ is a non-singular fan, then
$\XXX_{(\Sigma,\Sigma^0_{\extup})}$
is the toric variety $X_{\Sigma}$.
\end{enumerate}
\end{Remark}

\subsection{Torus Actions}
The torus action functor
\[
a_{(\Sigma,\Sigma^0)}:\XXX_{(\Sigma,\Sigma^0)} \times \Spec k[M] \longrightarrow \XXX_{(\Sigma,\Sigma^0)},
\]
is defined as follows.
Let $\xi=(\pi:\SSS \to \OO_X,\alpha:\MMM\to \OO_X,\eta:\SSS\to \MMM)$
be an object in $\XXX_{(\Sigma,\Sigma^0)}$.
Let $\phi:M \to \OO_X$ be a map of monoids
from a constant sheaf $M$ on $X$ to $\OO_X$, i.e., an $X$-valued point
of $\Spec k[M]$. Here $\OO_X$ is regarded as a sheaf of monoids
under multiplication.
We define $a_{(\Sigma,\Sigma^0)}(\xi,\phi)$ to be
$(\phi\cdot\pi:\SSS \to \OO_X,\alpha:\MMM\to \OO_X,\phi\cdot\eta:\SSS\to \MMM)$,
where $\phi\cdot\pi(s):=\phi(s)\cdot\pi(s)$ and
$\phi\cdot\eta(s):=\phi(s)\cdot\eta(s)$.
Let $h:\MMM_1\to \MMM_2$ be a morphism in $\XXX_{(\Sigma,\Sigma^0)}\times \Spec k[M]$
from $(\xi_1,\phi)$ to $(\xi_2,\phi)$, where
$\xi_i=(\pi:\SSS \to \OO_X,\alpha:\MMM_i\to \OO_X,\eta_i:\SSS\to \MMM_i)$
for $i=1,2$,
and $\phi:M \to \OO_X$ is an $X$-valued point
of $\Spec k[M]$. We define $a_{(\Sigma,\Sigma^0)}(h)$ to be $h$.
We remark that this action commutes with
the torus action of $\Spec k[M]$ on $X_{\Sigma}$.

\begin{Definition}
\label{torusequiv}
Let $(\Sigma_i\ \textup{in}\ N_{i,\RR},\Sigma_i^0)$ be a stacky fan
and $\XXX_{(\Sigma_i,\Sigma_i^0)}$ the associated toric algebraic stack
for $i=1,2$.
Put $M_i=\Hom_{\ZZ}(N_i,\ZZ)$.
Let us denote by $a_{(\Sigma_i,\Sigma_i^0)}:\XXX_{(\Sigma_i,\Sigma_i^0)}\times \Spec k[M_i]
\to \XXX_{(\Sigma_i,\Sigma_i^0)}$ the torus action.
A 1-morphism $f:\XXX_{(\Sigma_1,\Sigma_1^0)}\to \XXX_{(\Sigma_2,\Sigma_2^0)}$
is {\it torus-equivariant} if
the restriction $f_0$
of $f$ to $\Spec k[M_1]\subset \XXX_{(\Sigma_1,\Sigma_1^0)}$
defines a homomorphism of group $k$-schemes
$f_0:\Spec k[M_1]\to\Spec k[M_2]\subset \XXX_{(\Sigma_2,\Sigma_2^0)}$,
and the diagram
\[
\xymatrix{
 \XXX_{(\Sigma_1,\Sigma_1^0)}\times \Spec k[M_1]\ar[r]^{f\times f_0}\ar[d]^{a_{(\Sigma_1,\Sigma_1^0)}} & \XXX_{(\Sigma_2,\Sigma_2^0)}\times \Spec k[M_2] \ar[d]^{a_{(\Sigma_2,\Sigma_2^0)}} \\
 \XXX_{(\Sigma_1,\Sigma_1^0)} \ar[r]^{f} & \XXX_{(\Sigma_2,\Sigma_2^0)} \\
 }
\]
commutes in 2-categorical sense.
Similarly, we define the torus-equivariant (1-)morphisms
from a toric algebraic stack (or toric variety) to a toric algebraic stack (or toric variety).
We remark that $\pi_{(\Sigma,\Sigma^0)}:\XXX_{(\Sigma,\Sigma^0)} \to X_{\Sigma}$
is torus-equivariant.

\end{Definition}

\subsection{Algebraicity}

Now we recall some results which we use later.

\begin{Theorem}[(cf. \cite{I2} Theorem 4.5)]
\label{algebraic}
The stack $\XXX_{(\Sigma,\Sigma^0)}$ is a smooth Deligne-Mumford stack
of finite type and separated over $k$,
and the functor $\pi_{(\Sigma,\Sigma^0)}:\XXX_{(\Sigma,\Sigma^0)}\to X_{\Sigma}$ is a coarse moduli map.
\end{Theorem}

\begin{Remark}
\label{degefun}
By \cite{I2}, we can define the toric algebraic stack $\XXX_{(\Sigma,\Sigma^0)}$ over $\ZZ$. The stack $\XXX_{(\Sigma,\Sigma^0)}$ is an 
(not necessarily Deligne-Mumford) Artin stack over $\ZZ$.
In characteristic zero,
toric algebraic stacks are always Deligne-Mumford.

Let $f:\XXX_{(\Sigma,\Sigma^0)}\to \XXX_{(\Delta,\Delta^0)}$
be a functor (not necessarily torus-equivariant).
Then by the universality of coarse moduli spaces,
there exists a unique morphism
$f_c:X_{\Sigma}\to X_{\Delta}$
such that $f_c\circ \pi_{(\Sigma,\Sigma^0)}=
\pi_{(\Delta,\Delta^0)}\circ f$.
\end{Remark}

Here we shall give another presentation of $\XXX_{(\Sigma,\Sigma^0)}$,
that is
more directly represented in terms of logarithmic geometry.
It is important for the later proofs.
Let $(\UUU,\pi_{\UUU})$ be the universal pair on $X_{\Sigma}$
satisfying (1), (2), (3), (4) in Definition~\ref{torst},
which corresponds to
$\textup{Id}_{X_{\Sigma}}\in \Hom(X_{\Sigma},X_{\Sigma})$.
It follows from the construction in \cite{AMRT} that
the log structure associated to
$\pi_{\UUU}:\UUU\to \OOO_{X_{\Sigma}}$
is the canonical log structure $\MMM_{\Sigma}$ on $X_{\Sigma}$.
By \cite[4.4]{I2}, the stack $\XXX_{(\Sigma,\Sigma^0)}$ is naturally
isomorphic to
the stack $\XXX_{\Sigma}(\Sigma^0)$
over the toric variety $X_{\Sigma}$, which is defined as follows.
For any morphism $f:Y\to X_{\Sigma}$,
objects in $\XXX_{\Sigma}(\Sigma^0)$ over $f:Y\to X_{\Sigma}$
are morphisms of fine log schemes
$(f,\phi):(Y,\NNN) \to (X_{\Sigma},\MMM_{\Sigma})$
such that for every geometric point $\bar{y}\to Y$,
$\bar{\phi}:f^{-1}\bar{\UUU}_{\bar{y}}=f^{-1}\bar{\MMM}_{\Sigma,\bar{y}} \to \bar{\NNN}_{\bar{y}}$ is a $(\Sigma,\Sigma^0)$-free resolution.
(We shall call such a morphism
$(Y,\NNN)\to (X,\MMM_{\Sigma})$ a {\it $\Sigma^0$-FR morphism}.)
A morphism $(Y,\NNN)_{/(X_{\Sigma},\MMM_{\Sigma})}\to (Y',\NNN')_{/(X_{\Sigma},\MMM_{\Sigma})}$ in $\XXX_{\Sigma}(\Sigma^0)$ is a $(X_{\Sigma},\MMM_{\Sigma})$-morphism
$(\alpha,\phi):(Y,\NNN)\to (Y',\NNN')$ such that $\phi:\alpha^*\NNN'\to \NNN$ is an
isomorphism.

\begin{Remark}
Let $(f,h):(X,\MMM)\to (Y,\NNN)$ be a morphism of log schemes.
If
$h:f^*\NNN\to \MMM$ is an isomorphism,
we say that $(f,h)$ is {\it strict}.

\smallskip

We call $\XXX_{(\Sigma,\Sigma^0)}\cong\XXX_{\Sigma}(\Sigma^0)$ the toric algebraic stacks (or toric stacks) associated to $(\Sigma,\Sigma^0)$.
\end{Remark}

We shall collect some technical lemmata \ref{sub1}, \ref{sub2}, and
\ref{alg} (cf. \cite[2.16, 2.17, 3.4]{I2}),
which we will apply to the proof of Theorem~\ref{main} and \ref{main2}.
Let $(\Sigma,\Sigma^0)$
be a stacky fan. Assume that $\Sigma$ is a cone $\sigma$
such that $\dim \sigma=\rank N$, i.e., full-dimensional.
Set $P=\sigma^{\vee}\cap M$ ($M=\Hom_{\ZZ} (N,\ZZ)$).
The monoid $P$ is a simplicially toric sharp monoid,
and there is a natural isomorphism $X_{\Sigma}\cong\Spec k[P]$.

\begin{Lemma}
\label{sub1}
Let $r:P\to \NN^d$ be the minimal free resolution.
Let us denote by $e_{\rho}$ the irreducible element
in $\NN^d$ which corresponds to a ray $\rho$ in $\sigma$,
and let $t:\NN^d\to \NN^d$ be the map defined by $e_{\rho}\mapsto n_{\rho}\cdot e_{\rho}$ where $n_{\rho}$ is the level of $\sigma^0=\Sigma^0$ on $\rho$.
Let $(\Spec k[P],\MMM_P)$ and $(\Spec k[\NN^d],\MMM_{\NN^d})$ be
toric varieties with canonical log structures,
and let $(\pi,\eta):(\Spec k[\NN^d],\MMM_{\NN^d}) \to (\Spec k[P],\MMM_{P})$
the morphism of fine log schemes induced by $l:=t\circ r:P{\to}\NN^d{\to}\NN^d$. Then $(\pi,\eta)$ is a $\Sigma^0$-FR morphism.
\end{Lemma}

\begin{Lemma}
\label{sub2}
Let $(q,\gamma):(S,\NNN)\to (\Spec k[P],\MMM_{P})$ be a morphism of fine log schemes.
Let $c:P\to \MMM_P$ be a chart.
Let $\bar{s}$ be a geometric point on $S$.
Suppose that there exists
 a morphism $\xi:\NN^d\to \bar{\NNN}_{\bar{s}}$ such that
the composite $\xi\circ l:P\to \bar{\NNN}_{\bar{s}}$ is equal to
$\bar{\gamma}_{\bar{s}}\circ \bar{c}_{\bar{s}}:P\to q^{-1}\bar{\MMM}_{P,\bar{s}} \to \bar{\NNN}_{\bar{s}}$  $($with notation as in Lemma~\ref{sub1}$)$.
Assume that $\xi:\NN^d\to \bar{\NNN}_{\bar{s}}$ \'etale locally lifts to a
chart.
Then there exists an \'etale neighborhood $U$ of
$\bar{s}$ in which we have a chart $\varepsilon :\NN^d \to \NNN$ such that
the following diagram
\[
 \xymatrix{
 P \ar[d]_{c} \ar[r]^{l} & \NN^d \ar[d]^{\varepsilon}   \\
 q^*\MMM_{P} \ar[r]^{\gamma} & \NNN \\
  }
\]
commutes and the composite $\NN^d \stackrel{\varepsilon}{\to} \NNN \to \bar{\NNN}_{\bar{s}}$ is equal to $\xi$.
\end{Lemma}

Let $l:P\to \NN^d$ be the homomorphism as in Lemma~\ref{sub1}.
Let us denote by $G:=((\NN^d)^{gp}/l^{gp}(P^{gp}))^D$
the Cartier dual of the finite group $(\NN^d)^{gp}/l^{gp}(P^{gp})$.
The finite group scheme $G$ naturally acts on $\Spec k[\NN^d]$ as follows.
For a $k$-ring $A$, an $A$-valued point
$a:(\NN^d)^{gp}/l^{gp}(P^{gp}) \to A^*$ of $G$
sends an $A$-valued point $x:\NN^d\to A$ (a map of monoids)
of $\Spec k[\NN^d]$ to $a\cdot x:\NN^d\to A$; $n\mapsto a(n)\cdot x(n)$.
Since $G$ is \'etale over $k$ ($\textup{ch}(k)=0$),
the quotient stack $[\Spec k[\NN^d]/G]$ is a smooth Deligne-Mumford stack
(\cite[10.13]{LM}) whose coarse moduli space is $\Spec k[\NN^d]^{G}=\Spec k[P]$, where $k[\NN^d]^{G}\subset k[\NN^d]$ is the subring of functions
invariant under the action of $G$. The quotient $[\Spec k[\ZZ^d]/G]$
is an open representable substack of $[\Spec k[\NN^d]/G]$, which defines a torus embedding.

\begin{Proposition}
\label{alg}
There exists an isomorphism
$[\Spec k[\NN^d]/G] \to \XXX_{\Sigma}(\Sigma^0)$
of stacks over $\Spec k[P]$,
which  sends the torus in $[\Spec k[\NN^d]/G]$ onto that of $\XXX_{\Sigma}(\Sigma^0)$.
Moreover the natural composite $\Spec k[\NN^d] \to [\Spec k[\NN^d]/G] \to \XXX_{\Sigma}(\Sigma^0)$ corresponds to $(\Spec k[\NN^d],\MMM_{\NN^d})\to (\Spec k[P],\MMM_P)$.

\end{Proposition}

\subsection{Log structures on toric algebraic stacks}
Let $i_{(\Sigma,\Sigma^0)}:T_{\Sigma}=\Spec k[M]\to \XXX_{(\Sigma,\Sigma^0)}$
denote its torus embedding.
The complement $\DDD_{(\Sigma,\Sigma^0)}:=\XXX_{(\Sigma,\Sigma^0)}-T_{\Sigma}$
with reduced closed substack structure
is a normal crossing divisor (cf. \cite[4.17]{I2} or Proposition~\ref{alg}).
The stack $\XXX_{(\Sigma,\Sigma^0)}$
has the log structure
$\MMM_{(\Sigma,\Sigma^0)}$ arising from $\DDD_{(\Sigma,\Sigma^0)}$
 on the \'etale site
$\XXX_{(\Sigma,\Sigma^0),\textup{et}}$.
Moreover we have
\[
\MMM_{(\Sigma,\Sigma^0)}=\OOO_{\XXX_{(\Sigma,\Sigma^0)}}\cap i_{(\Sigma,\Sigma^0)*}\OOO_{T_{\Sigma}}^*\subset \OOO_{\XXX_{(\Sigma,\Sigma^0)}}
\]
where $\OOO_{\XXX_{(\Sigma,\Sigma^0)}}\cap i_{(\Sigma,\Sigma^0)*}\OOO_{T_{\Sigma}}^*$ denotes the subsheaf of $\OOO_{\XXX_{(\Sigma,\Sigma^0)}}$
consisting of regular functions on $\XXX_{(\Sigma,\Sigma^0)}$
whose restriction to $T_{\Sigma}$ is invertible.
The coarse moduli map $\pi_{(\Sigma,\Sigma^0)}:\XXX_{(\Sigma,\Sigma^0)}\to X_{\Sigma}$ induces a
morphism of log stacks,
$(\pi_{(\Sigma,\Sigma^0)},h_{(\Sigma,\Sigma^0)}):
(\XXX_{(\Sigma,\Sigma^0)},\MMM_{(\Sigma,\Sigma^0)})
\to (X_{\Sigma},\MMM_{\Sigma})$.
Here $h_{(\Sigma,\Sigma^0)}: \pi_{(\Sigma,\Sigma^0)}^*\MMM_{\Sigma}\to \MMM_{(\Sigma,\Sigma^0)}$ arises from the natural diagram
\[
\xymatrix{
\pi_{(\Sigma,\Sigma^0)}^{-1}\MMM_{\Sigma}\ar[d]\ar[r] & \MMM_{(\Sigma,\Sigma^0)}\ar[d]\\
 \pi_{(\Sigma,\Sigma^0)}^{-1}\OOO_{X_{\Sigma}}\ar[r] &\OOO_{\XXX_{(\Sigma,\Sigma^0)}}.\\
 }
\]
Similarly, a functor $f:\XXX_{(\Sigma,\Sigma^0)}\to
\XXX_{(\Delta,\Delta^0)}$
such that $f(T_{\Sigma})\subset T_{\Delta}\subset \XXX_{(\Delta,\Delta^0)}$
naturally induces the canonical homomorphism $h_f:f^*\MMM_{(\Delta,\Delta^0)} \to \MMM_{(\Sigma,\Sigma^0)}$ which is induced by
$f^{-1}\MMM_{(\Delta,\Delta^0)} \to f^{-1}\OOO_{\XXX_{(\Delta,\Delta^0)}}\to \OOO_{\XXX_{(\Sigma,\Sigma^0)}}$.
We shall refer to this homomorphism $h_{f}$ as
the {\it homomorphism induced by $f$}.

The log structure $\MMM_{(\Sigma,\Sigma^0)}$
on $\XXX_{\Sigma}(\Sigma^0)=\XXX_{(\Sigma,\Sigma^0)}$
has the following modular interpretation.
Let $f:Y\to \XXX_{\Sigma}(\Sigma^0)=\XXX_{(\Sigma,\Sigma^0)}$ be
a morphism from a $k$-scheme $Y$, which corresponds to a $\Sigma^0$-FR morphism
$(Y,\MMM_Y)\to (X_{\Sigma},\MMM_{\Sigma})$.
We attach the log structure $\MMM_Y$ to $f:Y\to \XXX_{\Sigma}(\Sigma^0)$,
and it gives rise to log structure $\MMM_{(\Sigma,\Sigma^0)}'$ on
$\XXX_{\Sigma}(\Sigma^0)=\XXX_{(\Sigma,\Sigma^0)}$.
We claim $\MMM_{(\Sigma,\Sigma^0)}=\MMM_{(\Sigma,\Sigma^0)}'$.
To see this, note the following observation.
Let $U\to \XXX_{\Sigma}(\Sigma^0)$ be an \'etale cover by a scheme $U$
and $\textup{pr}_1,\textup{pr}_2:U\times_{\XXX_{\Sigma}(\Sigma^0)}U\rightrightarrows U$
the \'etale groupoid.
A log structure on $\XXX_{\Sigma}(\Sigma^0)$ amounts to
a descent data $(\MMM_U,\textup{pr}_1^*\MMM_U\cong\textup{pr}_2^*\MMM_U)$ where $\MMM_U$ is a fine log structure on $U$.
Given a data $(\MMM_U,\textup{pr}_1^*\MMM_U\cong\textup{pr}_2^*\MMM_U)$, if $\MMM_U$ arises from a normal crossing divisor on $U$,
then $\MMM_U\subset \OOO_U$ and $\textup{pr}_1^*\MMM_U=\textup{pr}_2^*\MMM_U\subset \OOO_{U\times_{\XXX_{\Sigma}(\Sigma^0)}U}$.
By Proposition~\ref{alg}, there is an \'etale cover 
$f:U\to \XXX_{\Sigma}(\Sigma^0)$ such that $f^*\MMM_{(\Sigma,\Sigma^0)}'$ arises from
the divisor $f^{-1}(\DDD_{(\Sigma,\Sigma^0)})$. Then from the above observation
and  the equality $f^*\MMM_{(\Sigma,\Sigma^0)}=f^*\MMM_{(\Sigma,\Sigma^0)}'\subset \OOO_U$, we conclude that
$\MMM_{(\Sigma,\Sigma^0)}'$ is isomorphic to $\MMM_{(\Sigma,\Sigma^0)}$
up to unique isomorphism.
For generalities concerning log structures on stacks,
we refer to (\cite[Section 5]{OL}).

\begin{Remark}
The notion of stacky fans was introduced
in \cite[Section 3]{BCS}. 
We should remark that in \cite{BCS}, given a stacky fan $(\Sigma,\Sigma^0)$
whose rays in $\Sigma$ span the vector space
$N_{\RR}$, Borisov-Chen-Smith constructed a smooth Deligne-Mumford stack 
over $\CC$ whose coarse
moduli space is the toric variety $X_{\Sigma}$,
called the toric Deligne-Mumford stack.
Their approach is a generalization of the
global quotient constructions of toric varieties due to D. Cox.
However, it seems quite difficult to show
the 2-category (or the associated 1-category) of toric Deligne-Mumford stacks in the sense of \cite{BCS}
is equivalent to the category of stacky fans
by their machinery.
In Section 5, we explain the relationship with \cite{BCS}.
\end{Remark}

\section{The proof of Theorem~\ref{main}}
In this Section, we shall prove Theorem~\ref{main}.
As in Section 2, we continue to work over the fixed base field $k$ of characteristic zero.
The proof proceeds in several steps.

\begin{Lemma}
\label{logmap}
Let $\XXX_{(\Sigma,\Sigma^0)}$ and $\XXX_{(\Delta,\Delta^0)}$ be
toric algebraic stacks arising from stacky fans $(\Sigma\ \textup{in}\ N_{1,\RR},\Sigma^0)$
and $(\Delta\ \textup{in}\ N_{2,\RR},\Delta^0)$ respectively.
Let $f:\XXX_{(\Sigma,\Sigma^0)}\to \XXX_{(\Delta,\Delta^0)}$
be a functor such that $f(T_{\Sigma})\subset T_{\Delta}\subset \XXX_{(\Delta,\Delta^0)}$.
Let $f_c:X_{\Sigma}\to X_{\Delta}$ be the morphism induced by $f$
$($cf. Remark~\ref{degefun}$)$.
Then there exists a natural commutative diagram
of log stacks
\[
\begin{CD}
(\XXX_{(\Sigma,\Sigma^0)},\MMM_{(\Sigma,\Sigma^0)}) @>{(f,h_f)}>> (\XXX_{(\Delta,\Delta^0)},\MMM_{(\Delta,\Delta^0)}) \\
@V{(\pi_{(\Sigma,\Sigma^0)},h_{(\Sigma,\Sigma^0)})}VV  @VV{(\pi_{(\Delta,\Delta^0)},h_{(\Delta,\Delta^0)})}V \\
(X_{\Sigma},\MMM_{\Sigma}) @>{(f_c,h_{f_c})}>> (X_{\Delta},\MMM_{\Delta}). \\
\end{CD}
\]
\end{Lemma}

\Proof
We use the same notation as in Section 2.4.
Note that $f_c$
commutes with torus embeddings.
We define a homomorphism
$h_{f_c}:f_c^*\MMM_{\Delta} \to
\MMM_{\Sigma}$
to be the homomorphism induced by
$f_c^{-1}\MMM_{\Delta}\to f_c^{-1}\OOO_{X_{\Delta}}\to \OOO_{X_{\Sigma}}$.
Then since $h_{f}$, $h_{f_c}$, $h_{(\Sigma,\Sigma^0)}$, and
$h_{(\Delta,\Delta^0)}$ are induced by the homomorphisms of structure
sheaves (Section 2.4),
$h_{(\Sigma,\Sigma^0)}\circ \pi_{(\Sigma,\Sigma^0)}^*h_{f_c}:(f_c\circ \pi_{(\Sigma,\Sigma^0)})^*\MMM_{\Delta}\to \MMM_{(\Sigma,\Sigma^0)}$
is equal to
$h_{f}\circ f^*h_{(\Delta,\Delta^0)}:(\pi_{(\Delta,\Delta^0)}\circ f)^*\MMM_{\Delta}\to \MMM_{(\Sigma,\Sigma^0)}$.
Thus we have the desired diagram.
\QED

\begin{Proposition}
\label{eqeq}
With notation in Lemma~\ref{logmap}, 
if $f$ is torus-equivariant, then the morphism $f_c$
is torus-equivariant. Moreover the morphism $f_c$ corresponds to
the map of fans $L:\Sigma \to \Delta$ such that $L(\Sigma^0) \subset \Delta^0$.
\end{Proposition}

\Proof
Clearly, the restriction of $f_c$ to $T_{\Sigma}$ induces
a homomorphism of group $k$-schemes
$T_{\Sigma}\to T_{\Delta}\subset X_{\Delta}$.
Note that $\XXX_{(\Sigma,\Sigma^0)}\times T_{\Sigma} \stackrel{f\times (f|_{T_{\Sigma}})}{\longrightarrow} \XXX_{(\Delta,\Delta^0)}\times T_{\Delta}\stackrel{a_{(\Delta,\Delta^0)}}{\longrightarrow} \XXX_{(\Delta,\Delta^0)}$
is isomorphic to $\XXX_{(\Sigma,\Sigma^0)}\times T_{\Sigma}\stackrel{a_{(\Sigma,\Sigma^0)}}{\longrightarrow} \XXX_{(\Sigma,\Sigma^0)} \stackrel{f}{\to} \XXX_{(\Delta,\Delta^0)}$.
Since $X_{(\Sigma,\Sigma^0)}\times T_{\Sigma}$ 
(resp. $X_{(\Delta,\Delta^0)}\times T_{\Delta}$)
is a coarse moduli space for $\XXX_{(\Sigma,\Sigma^0)}\times T_{\Sigma}$
(resp. $\XXX_{(\Delta,\Delta^0)}\times T_{\Delta}$),
thus $f_c$ is torus-equivariant.
Set $M_i=\Hom_{\ZZ} (N_i,\ZZ)$ for $i=1,2$.
Let $L^{\vee}:M_2\to M_1$ be the homomorphism of abelian groups
that is induced by the homomorphism of
group $k$-schemes $f_c|_{T_{\Sigma}}:T_{\Sigma}\to T_{\Delta}$.
The dual map $L:N_1=\Hom_{\ZZ} (M_1,\ZZ) \to \Hom_{\ZZ} (M_2,\ZZ)=N_2$
yields the map of fans $L_{\RR}:\Sigma\ \textup{in}\ N_{1,\RR} \to \Delta\ \textup{in}\ N_{2,\RR}$, which corresponds to the morphism $f_c$.
To complete the proof of this Proposition, it suffices to
show $L(\Sigma^0) \subset \Delta^0$.
To do this, we may assume $(\Sigma,\Sigma^0)=(\sigma,\sigma^0)$ and $(\Delta,\Delta^0)
=(\delta,\delta^0)$ where $\sigma$ and $\delta$ are cones.
We need the following lemma.
\QED

\begin{Lemma}
\label{constfreeres}
If $\sigma$ is a full-dimensional cone, then the $(\sigma,\sigma^{0})$-free resolution $($cf. Definition~\ref{torst}$)$ of $\sigma^{\vee}\cap M_1$
is given by
\[
\sigma^{\vee}\cap M_1\to \{ m \in M_1\otimes_{\ZZ}\QQ \|\ \langle m,n\rangle \in \ZZ_{\ge 0}\ \textup{for\ any}\ n \in \sigma^0 \}.
\]
\end{Lemma}

{\it Proof of Lemma.}
Let $P:=\sigma^{\vee}\cap M_1$.
We first show the case of $\sigma^{0}=\sigma_{\extup}^0$. We assume $\sigma^{0}=\sigma_{\extup}^0$.
Let $d$ be the rank of $M_1$.
Here $M_1=P^{gp}$.
Let $S$ be a submonoid in $\sigma\cap N_1$ which is generated by
the first lattice points of rays in $\sigma$.
Namely, $S=\sigma_{\extup}^0$ and $S\cong\NN^d$.
Put $F:=\{ h \in M_1\otimes_\ZZ \QQ \mid \langle h,s \rangle \in \ZZ_{\ge 0}$ for any $s \in S\}$.
It is clear that
$F$ is isomorphic to $\NN^d$.
 Since $\sigma^{\vee}\cap M_1=\{ h \in M_1\otimes_\ZZ \QQ \mid \langle h,s\rangle \in \ZZ_{\ge 0}$ for any $s \in \sigma\cap N_1\}$, we have
$P \subset F \subset M_1\otimes_{\ZZ}\QQ$.
We will show that the natural injective map $i:P\to F$
is the minimal free resolution (cf. Definition~\ref{mini}).
Since the monoid $\sigma \cap N_1$ is a fine sharp monoid,
it has only finitely many irreducible elements and it is
generated by the irreducible elements (\cite[Lemma 3.9]{OL}).
Let $\{ \zeta_1,\ldots,\zeta_r\}$ (resp. $\{ e_1,\ldots,e_d\}$) be irreducible elements of $\sigma\cap N_1$ (resp. $F$).
For each irreducible element $e_i$ in $F$,
we put $\langle e_i, \zeta_j\rangle=a_{ij}/b_{ij} \in \QQ$ with some $a_{ij}\in \ZZ_{\ge 0}$
and $b_{ij} \in \NN$.
Then we have
\[
(\Pi_{0\le j \le r}b_{ij})\cdot e_i(\sigma\cap N_1) \subset \ZZ_{\ge0}
\]
and thus $i:P \to F$ satisfies the property (1) of Definition~\ref{mini}.
To show our claim, it suffices to prove that the $i:P\to F$
satisfies the property (2) of Definition~\ref{mini}.
Let $j:P\to G$ be an injective homomorphism of monoids such that
$j(P)$ is close to $G$ and $G\cong \NN^d$.
The monoid $P$ has the natural injection $l:P\to M_1\otimes_{\ZZ}{\QQ}$.
On the other hand, for any element $e$ in $G$, there exists a positive integer
$n$ such that $n\cdot e \in j(P)$. Therefore we have a unique homomorphism
 $\alpha:G\to M_1\otimes_{\ZZ}\QQ$ which extends $l:P\to M_1\otimes_{\ZZ}{\QQ}$ to
$G$.
We claim that there exists a sequence of inclusions
\[
P \subset F \subset \alpha(G) \subset M_1\otimes_{\ZZ}\QQ.
\]
If we put $\alpha(G)^{\vee}:=\{f \in N_1\otimes_{\ZZ}{\QQ}\mid \langle p,f\rangle\in\ZZ_{\ge0}$ for any $p \in \alpha(G) \} \cong \NN^{d}$
and
$F^{\vee}:=\{f \in N_1\otimes_{\ZZ}{\QQ}\mid \langle p,f\rangle\in\ZZ_{\ge0}$ for any $p \in F \}=S \cong \NN^{d}$, then our claim is equivalent to
the claim $\alpha(G)^{\vee} \subset F^{\vee}$.
However the latter claim is clear. Indeed, the sublattice
 $S$ is the maximal sublattice
of $\sigma\cap N_1$ which is free and 
close to $\sigma\cap N_1$.
The sublattice
$\alpha(G)^{\vee}$ is also close to $\sigma\cap N_1$, and thus each irreducible generator of $\alpha(G)^{\vee}$ lies on a ray of $\sigma$.

Finally, we consider the general case.
Put $F_{\extup}:=\{ h \in M_1\otimes_\ZZ \QQ \mid \langle h,s \rangle \in \ZZ_{\ge 0}$ for any $s \in \sigma_{\extup}^0\}$
and
$F:=\{ h \in M_1\otimes_\ZZ \QQ \mid \langle h,s \rangle \in \ZZ_{\ge 0}$ for any $s \in \sigma^0\}$. (Note that notation changed.) Then we have a natural
injective map $F_{\extup}\to F$ because $\sigma^0\subset \sigma^0_{\extup}$. Given a ray $\rho\in \sigma(1)$,
the corresponding irreducible element of $F_{\extup}$
(resp. $F$) (cf. Lemma~\ref{irr}) is $m_{\rho}\ (\textup{resp.}\ m'_{\rho}) \in M_1\otimes_{\ZZ}\QQ$
such that $\langle m_{\rho},\zeta_{\rho}\rangle=1$ 
(resp. $\langle m'_{\rho},n_{\rho}\zeta_{\rho}\rangle=1$),
and $\langle m_{\rho},\zeta_{\xi}\rangle=0$ (resp. $\langle m'_{\rho},\zeta_{\xi}\rangle=0$) for any ray $\xi$ with $\xi\neq\rho$.
Here for each ray $\alpha$, $\zeta_{\alpha}$ denotes the first lattice point of $\alpha$,
and $n_{\alpha}$ denotes the level of $\sigma^0$ on $\alpha$.
The natural injection $F_{\extup}\to F$
identifies $m_{\rho}$ with $n_{\rho}m'_{\rho}$. Thus $\sigma^{\vee}\cap M_1\to F$ is a $(\sigma,\sigma^0)$-free resolution of $\sigma^{\vee}\cap M_1$.
\QED

We continue the proof of Proposition.
We shall assume $L(\Sigma^0) \nsubseteq \Delta^0$
and show that such an assumption gives rise to a contradiction.
First we show a contradiction for the case when
$\sigma$ and $\delta$ are full-dimensional cones.
Set $P:=\sigma^{\vee}\cap M_1$ and $Q:=\delta^{\vee}\cap M_2$.
Note that since $\sigma$ and $\delta$ are full-dimensional,
$P$ and $Q$ are sharp (i.e., unit-free).
Let us denote by $o$ (resp. $o'$) the origin of $\Spec k[P]$ 
(which corresponds to the ideal $(P)$) (resp. the origin of $\Spec k[Q]$).
Then $f_c$ sends $o$ to $o'$.
Consider the composite
$\alpha:\Spec \OOO_{\Spec k[\NN^d],\bar{s}} \to \Spec k[\NN^d]\to
[\Spec k[\NN^d]/G]\cong\XXX_{(\Sigma,\Sigma^0)}$
of natural morphisms (cf. Proposition~\ref{alg}),
where $s$ is the origin of $\Spec k [\NN^d]$.
Then by Lemma~\ref{logmap},
there exists the following commutative diagram
\[
\xymatrix{
 M_2 \ar[d]_{L^{\vee}} & Q=\alpha^{-1}\pi_{(\Sigma,\Sigma^0)}^{-1}f_c^{-1}\bar{\MMM}_{\Delta} \ar[l] \ar[r] \ar[d] & \alpha^{-1}f^{-1}\bar{\MMM}_{(\Delta,\Delta^0)} \ar[d]  \\M_1  & P=\alpha^{-1}\pi_{(\Sigma,\Sigma^0)}^{-1}\bar{\MMM}_{\Sigma} \ar[l] \ar[r] & \alpha^{-1}\bar{\MMM}_{(\Sigma,\Sigma^0)}.   \\
  }
\]
On the other hand,
set $F:=\{ m \in M_1\otimes_{\ZZ}\QQ \|\ \langle m,n\rangle \in \ZZ_{\ge 0} \textup{for\ any} n \in \sigma^0 \}$
and 
$F':=\{ m \in M_2\otimes_{\ZZ}\QQ \|\ \langle m,n\rangle \in \ZZ_{\ge 0} \textup{for\ any} n \in \delta^0 \}$.
Then by the above Lemma, the $(\sigma,\sigma^{0})$-free resolution
$\alpha^{-1}\pi_{(\Sigma,\Sigma^0)}^{-1}\bar{\MMM}_{\Sigma}\to
\alpha^{-1}\bar{\MMM}_{(\Sigma,\Sigma^0)}$ is identified with
the natural inclusion $P:=\sigma^{\vee}\cap M_1 \hookrightarrow F$
(the monoid $\alpha^{-1}\bar{\MMM}_{(\Sigma,\Sigma^0)}$ can be canonically embedded
into $M_1\otimes\QQ$).
Similarly, $\alpha^{-1}\pi_{(\Sigma,\Sigma^0)}^{-1}f_c^{-1}\bar{\MMM}_{\Delta} \to \alpha^{-1}f^{-1}\bar{\MMM}_{(\Delta,\Delta^0)}$
can be identified with the natural inclusion $Q=\delta^{\vee}\cap M_2 \hookrightarrow F'$.
The homomorphism $\alpha^{-1}f^{-1}\bar{\MMM}_{(\Delta,\Delta^0)}\to\alpha^{-1}\bar{\MMM}_{(\Sigma,\Sigma^0)}$
can be naturally embedded into $L^{\vee}\otimes\QQ:M_2\otimes\QQ\to M_1\otimes\QQ$.
However the assumption $L(\Sigma^0) \nsubseteq \Delta^0$
implies $L^{\vee}(F') \nsubseteq  F$. It gives rise to a contradiction.
Next consider the general case, i.e., $\sigma$ and $\delta$
are not necessarily full-dimensional.
Choose splittings $N_i\cong N_i'\oplus N_i''$ (i=1,2),
$\sigma\cong \sigma'\oplus \{ 0\}$, $\delta\cong \delta'\oplus \{ 0 \}$
such that $\sigma'$ and $\delta'$ are full-dimensional in $N_{1,\RR}'$ and
$N_{2,\RR}'$, respectively.
Note that
$\XXX_{(\sigma,\sigma^0)}\cong \XXX_{(\sigma',\sigma^{'0})}\times \Spec k[M_1'']$ and
$\XXX_{(\delta,\delta^0)}\cong \XXX_{(\delta',\delta^{'0})}\times \Spec k[M_2'']$.
Consider the following sequence of torus-equivariant morphisms
\[
\XXX_{(\sigma',\sigma^{'0})}\stackrel{i}{\to}\XXX_{(\sigma,\sigma^0)}\stackrel{f}{\to}
\XXX_{(\delta,\delta^0)}\cong\XXX_{(\delta',\delta^{'0})}\times \Spec k[M_2''] \stackrel{\textup{pr}_1}{\to} \XXX_{(\delta',\delta^{'0})},
\]
where $i$ is determined by the natural inclusion
$N_1'\hookrightarrow N_1'\oplus N_1''$,
and $\textup{pr}_1$ is the first projection.
Notice that $i$ and $\textup{pr}_1$ naturally induce
isomorphisms $i^*\MMM_{(\sigma,\sigma^0)}\stackrel{\sim}{\to}\MMM_{(\sigma',\sigma^{'0})}$ and $\textup{pr}_1^*\MMM_{(\delta',\delta^{'0})}\stackrel{\sim}{\to} \MMM_{(\delta,\delta^0)}$ respectively.
Thus the general case follows from the full-dimensional case.
\QED

\begin{Remark}
\label{des}
(1) By Proposition~\ref{eqeq},
there exists the natural functor
\[
c:\mathfrak{Torst}\to \textup{Simtoric},\ \XXX_{(\Sigma,\Sigma^0)}\mapsto
X_{\Sigma}
\]
which sends a torus-equivariant morphism $f:\XXX_{(\Sigma,\Sigma^0)}\to \XXX_{(\Delta,\Delta^0)}$
to $f_c:X_{\Sigma}\to X_{\Delta}$, where $f_c$
is the unique morphism such that $f_c\circ \pi_{(\Sigma,\Sigma^0)}=
\pi_{(\Delta,\Delta^0)}\circ f$.

(2) We can define a 2-functor 
\[
\Phi:\mathfrak{Torst}\to
(\textup{Category\ of\ stacky\ fans}),\ \XXX_{(\Sigma,\Sigma^0)}\mapsto(\Sigma,\Sigma^0)
\]
as follows.
For each torus-equivariant morphism
$f:\XXX_{(\Sigma,\Sigma^0)}\to \XXX_{(\Delta,\Delta^0)}$,
the restriction of $f$ to the torus $\Spec k[M_1]\subset \XXX_{(\Sigma,\Sigma^0)}$
induces the homomorphism $\Phi(f):N_1\to N_2$
that defines a morphism of stacky fans
$\Phi(f):(\Sigma,\Sigma^0)\to (\Delta,\Delta^0)$ by Proposition~\ref{eqeq}.
For each 2-isomorphism morphism $g:f_1\to f_2$
($f_1, f_2:\XXX_{(\Sigma,\Sigma^0)}\rightrightarrows \XXX_{(\Delta,\Delta^0)}$
are torus-equivariant morphisms),
define $\Phi(g)$ to be $\textup{Id}_{\Phi(f_1)}$. (Note that $\Phi(f_1)=\Phi(f_2)$.)
\end{Remark}

In order to show Theorem~\ref{main}, we show the following key Proposition.

\begin{Proposition}
\label{const}
Let $\xi:(\Sigma\ \textup{in}\ N_{1,\RR},\Sigma^0)\to (\Delta\ \textup{in}\ N_{2,\RR},\Delta^0)$ be
a morphism of stacky fans.
Let $(f,h_f):(X_{\Sigma},\MMM_{\Sigma})\to (X_{\Delta},\MMM_{\Delta})$ be 
the morphism of log toric varieties
induced by $\xi:\Sigma\to\Delta$.
Let $S$ be a $k$-scheme
and let $(\alpha,h):(S,\NNN)\to (X_{\Sigma},\MMM_{\Sigma})$
be a $\Sigma^0$-FR morphism.
Then there exist a fine log structure $\AAA$
on $S$, and morphisms of log structures
$a:\alpha^*f^*\MMM_{\Delta} \to \AAA$
and $\theta:\AAA\to \NNN$ which make the following diagram
\[
\begin{CD}
\alpha^*f^*\MMM_{\Delta} @>{a}>> \AAA \\
@V{\alpha^*h_f}VV  @V{\theta}VV \\
\alpha^*\MMM_{\Sigma}@>{h}>> \NNN, \\
\end{CD}
\]
commutative
and make $(f\circ \alpha,a):(S,\AAA)\to (X_{\Delta},\MMM_{\Delta})$
a $\Delta^0$-FR morphism.
The triple $(\AAA,a,\theta)$ is unique in the following sense:
If there exists such another triple $(\AAA',a',\theta')$,
then there exists a unique isomorphism $\eta:\AAA\to \AAA'$
which makes the diagram $(\clubsuit)$
\[
\xymatrix@R=3mm @C=12mm{
 \alpha^*f^*\MMM_{\Delta} \ar[dd]^{\alpha^*h_f} \ar[rr]^{a} \ar[rd]^{a'} &   & \AAA \ar[ld]^{\eta} \ar[dd]^{\theta}\\
  & \AAA' \ar[dr]^{\theta'}  &  \\
 \alpha^*\MMM_{\Sigma} \ar[rr]^{h} &   & \NNN \\
}
\]
commutative.
\end{Proposition}
We first show our claim for the case 
of $S=\Spec R$ where $R$ is a strictly Henselian local $k$-ring.
Note that if $\MMM$ is a fine saturated log structure on $S=\Spec R$,
then by \cite[Proposition 2.1]{OL},
there exists a chart $\bar{\MMM}(S)\to \MMM$ on $S$.
The chart induces an isomorphism
$\bar{\MMM}(S)\oplus R^*\stackrel{\sim}{\to}\MMM(S)$.
If a chart of $\MMM$ is fixed,
we abuse notation and usually write $\bar{\MMM}(S)\oplus R^*$
for the log structure $\MMM$.
Similarly, we write simply $\bar{\MMM}$ for $\bar{\MMM}(S)$.
Before the proof of Proposition, we prove the following lemma.

\begin{Lemma}
\label{monoid}
Set $P:=\alpha^{-1}\bar{\MMM}_{\Sigma}(S)$,
$Q:=\alpha^{-1}f^{-1}\bar{\MMM}_{\Delta}(S)$
and 
$\iota:=\alpha^{-1}\bar{h}_f:Q\to P$.
Let $\gamma:Q\stackrel{r}{\to}\NN^r \stackrel{t}{\to} \NN^r$ be the composite
map where $r$ is the minimal free resolution and
$t$ is the map defined as follows. For the irreducible element
$e_{\rho}\in \NN^r$ which corresponds to each ray $\rho$ in $\Delta$,
$t$ sends $e_{\rho}$ to $n_{\rho}\cdot e_{\rho}$ where
$n_{\rho}$ is the level of $\Delta^0$ on $\rho$.
Then there exists a unique homomorphism of monoids $l:\NN^r\to \bar{\NNN}$
such that the following diagram
\[
\begin{CD}
Q @>{\gamma}>> \NN^r \\
@V{\iota}VV   @V{l}VV \\
P @>{\bar{h}}>> \bar{\NNN} \\
\end{CD}
\]
commutes.
\end{Lemma}

\Proof
The uniqueness of $l$ follows from
the facts that $\bar{\NNN}$ is free
and
$\gamma(Q)$ is close to $\NN^r$.
To show the existence of $l$,
we may assume that $\Sigma$ and $\Delta$ are cones.
Set $\sigma=\Sigma,\ \sigma^0=\Sigma^0,\ \delta=\Delta$, and
$\delta^0=\Delta^0$.
Choose splittings 
$N_i\cong N_i'\oplus N_i''$ (i=1,2),
$\sigma\cong \sigma'\oplus \{ 0\}$, and
$\delta\cong \delta'\oplus \{ 0 \}$
such that $\dim N_1'=\dim \sigma'$ and $\dim N_2'=\dim\delta'$.
Then the projections $N_i'\oplus N_i''\to N_i'\ (i=1,2)$
yield the commutative diagram of log schemes
\[
\begin{CD}
(X_{\sigma},\MMM_{\sigma}) @>{(f,h_f)}>> (X_{\delta},\MMM_{\delta}) \\
@VVV   @VVV \\
(X_{\sigma'},\MMM_{\sigma'}) @>{(g,h_g)}>> (X_{\delta'},\MMM_{\delta'}), \\
\end{CD}
\]
where vertical arrows are strict morphisms induced by projections, and
$(g,h_g)$ is the morphism induced by $\xi |_{N_1'}:N_1'\to N_2'$.
Thus we may assume that $\sigma$ and $\delta$ are full-dimensional in $N_{1,\RR}$ and $N_{2,\RR}$ respectively.
Then
$P':=\sigma^{\vee}\cap M_1$ and $Q':=\delta^{\vee}\cap M_2$
are sharp (i.e. unit-free).
Let $R_1:P'\to \NN^{m}$ (resp. $R_2:Q'\to \NN^{n}$) be the
$(\sigma,\sigma^{0})$-free 
(resp. $(\delta,\delta^{0})$-free) resolution.
Let $f^{\#}:Q'\to P'$ be the homomorphism arising from $f$.
Then by the assumption $\xi(\Sigma^0)\subset \Delta^0$,
there exists a homomorphism $w:\NN^{n}\to \NN^{m}$
such that $w\circ R_2=R_1\circ f^{\#}$.
Taking Lemma~\ref{sub1} into account,
we see that our claim follows.\QED

{\it Proof of Proposition~\ref{const}}.
By Lemma~\ref{monoid}, there exists
a unique homomorphism $l:\NN^r\to \bar{\NNN}$.
By \cite[Proposition 2.1]{OL},
there exists a chart
$c_{\bar{\NNN}}:\bar{\NNN}\to \NNN$.
Then
maps $c_{\bar{\NNN}}\circ l$ and $c_{\bar{\NNN}}\circ l\circ \gamma$ induce the log structures
$\NN^r\oplus R^*$ and $Q\oplus R^*$ respectively.
On the other hand, by \cite[Proposition 2.1]{OL}
there exists a chart $c'_Q:Q\to \alpha^*f^*\MMM_{\Delta}$.
Let $i:R^*\hookrightarrow \NNN$ be the canonical
immersion. Let us denote by $j$ the composite map
\[
Q\stackrel{c'_Q}{\to}  \alpha^*f^*\MMM_{\Delta} \stackrel{\alpha^*h_f}{\to}\MMM_{\Sigma} \stackrel{h}{\to} \NNN \stackrel{(c_{\bar{\NNN}}\oplus i)^{-1}}{\to}\bar{\NNN}\oplus R^*\stackrel{\textup{pr}_2}{\to} R^*,
\]
and
define a chart $c_Q:Q\to \alpha^*f^*\MMM_{\Delta}$
by $Q\ni q\mapsto c_Q'(q)\cdot j(q)^{-1} \in\alpha^*f^*\MMM_{\Delta}$
(here $j(q)$ is viewed as an element in $\alpha^*f^*\MMM_{\Delta}$).
 Then we have the following commutative
diagram
\[
\xymatrix@R=3mm @C=6mm{
  & &  & \bar{\NNN}\oplus R^* \ar[ddd]_(0.3){c_{\bar{\NNN}}\oplus i}\ar@{=}[rr]&  & \bar{\NNN}\oplus R^* \ar[ddd]\\
  &  & &  & &  \\
  Q\oplus R^*\ar@{=}[rr]\ar[ddd]_{c_Q\oplus i_{\Delta}}&  & Q\oplus R^* \ar[rr]_(0.3){\gamma\oplus \textup{Id}_{R^*}}\ar[ddd]&  & \NN^r\oplus R^* \ar[ddd]\ar[uur]_{l\oplus \textup{Id}_{R^*}} &  \\
  & \alpha^*\MMM_{\Sigma} \ar[rr]^(0.3){h}&  & \NNN \ar[rr]& & \bar{\NNN} \\
  & & & & & \\
 \alpha^*f^*\MMM_{\Delta} \ar[rr] \ar[uur]^{\alpha^*h_f}&  & Q \ar[rr]^{\gamma}&  & \NN^r\ar[uur]^{l} &  \\
}
\]
where $i_{\Delta}:R^*\hookrightarrow \alpha^*f^*\MMM_{\Delta}$ is the canonical immersion.
The map $c_Q\oplus i_{\Delta}$ is an isomorphism, and ($\gamma\oplus \textup{Id}_{R^*})\circ (c_Q\oplus i_{\Delta})^{-1}:\alpha^*f^*\MMM_{\Delta}\to \NN^r\oplus R^*$
makes $f\circ \alpha:S\to X_{\Delta}$
a $\Delta^0$-FR morphism.
Thus
we have the desired diagram. 
Next we shall prove the uniqueness.
To prove this, as above, we fix the chart
$c_{\bar{\NNN}}:\bar{\NNN} \to \NNN$.
Suppose that  for $\lambda=1,2$, there exist a $\Delta^0$-FR morphism $a_\lambda:\alpha^*f^*\MMM_{\Delta}\to \AAA_{\lambda}$,
and a morphism of log structures $\theta_{\lambda}:\AAA_{\lambda}\to \NNN$ such that
$h\circ \alpha^*h_f=\theta_{\lambda}\circ a_{\lambda}$.
By the above argument (for the proof of the existence), we have a chart $c_Q:Q\to \alpha^*f^*\MMM_{\Delta}$ such that the image of the composite
\[
Q\stackrel{c_Q}{\to}  \alpha^*f^*\MMM_{\Delta} \stackrel{\alpha^*h_f}{\to}\alpha^*\MMM_{\Sigma} \stackrel{h}{\to} \NNN \stackrel{(c_{\bar{\NNN}}\oplus i)^{-1}}{\to}\bar{\NNN}\oplus R^*\stackrel{\textup{pr}_2}{\to} R^*
\]
is trivial. By Lemma~\ref{sub2}, $c_Q$ can be extended to a
chart $c_{\lambda}:\NN^r=\bar{\AAA}_{\lambda}\to \AAA_{\lambda}$
such that $c_{\lambda}\circ \gamma=a_{\lambda}\circ c_Q$ for $\lambda=1,2$.
Then the composite $\textup{pr}_2\circ (c_{\bar{\NNN}}\oplus i)^{-1}\circ \theta_{\lambda}\circ c_{\lambda}:\NN^r\to R^*$
induces a character $\textup{ch}_{\lambda}:(\NN^r)^{gp}/\gamma(Q)^{gp}\to R^*$.
Note that if $i_{\lambda}$
denotes the canonical immersion $R^*\hookrightarrow \AAA_{\lambda}$
for $\lambda=1,2$,
then $c_{\lambda}\oplus i_{\lambda}:\NN^r\oplus R^*\rightarrow \AAA_{\lambda}$
is an isomorphism.
Let us denote by $\eta:\AAA_1\to \AAA_2$ an isomorphism of log structures,
defined to
be the composite $\AAA_1\stackrel{(c_1\oplus i_1)^{-1}}{\to} \NN^r\oplus R^*\stackrel{\omega}{\to} \NN^r\oplus R^*\stackrel{(c_2\oplus i_2)}{\to} \AAA_2$,
where $\omega:\NN^r\oplus R^*\ni (n,u) \mapsto (n,\textup{ch}_{1}(n)\cdot \textup{ch}_{2}(n)^{-1}\cdot u) \in \NN^r\oplus R^*$.
Then it is easy to see that $\eta :\AAA_1\to \AAA_2$
 is a unique isomorphism
which makes all diagrams commutative.\QED

Next consider the case of a general $k$-scheme $S$.
First, we shall prove the uniqueness part.
If there exists a diagram as in $(\clubsuit)$ but without $\eta$, then
by the case of the spectrum of strictly Henselian local $k$-ring,
for every geometric point $\bar{s}$ on $S$,
there exists a unique homomorphism
$\eta_{\bar{s}}:\AAA_{\bar{s}}\to \AAA'_{\bar{s}}$
which makes the diagram $(\clubsuit)$
over $\Spec \OOO_{S,\bar{s}}$
commutative.
Thus, to prove the uniqueness, it suffices to show that
$\eta_{\bar{s}}$ can be extended to an isomorphism 
on some \'etale neighborhood
of $\bar{s}$,
which makes
the diagram $(\clubsuit)$ commutative.
To this aim, put $Q=\alpha^{-1}f^{-1}\bar{\MMM}_{\Delta,\bar{s}}$,
and choose a chart $c_Q:Q\to \alpha^{*}f^{*}{\MMM}_{\Delta}$
on some \'etale neighborhood $U$ of $\bar{s}$ (such an
existence follows from \cite[Proposition 2.1]{OL}).
We view the monoid $Q$ as a submonoid of $\NN^r\cong\bar{\AAA}_{\bar{s}}\cong \bar{\AAA}'_{\bar{s}}$.
Taking Lemma~\ref{sub2} and the existence of
$\eta_{\bar{s}}$ into account, after shrinking $U$ if necessary,
we can choose charts
$\NN^r \stackrel{c}{\to} \AAA$ and $\NN^r \stackrel{c'}{\to} \AAA'$ on $U$,
such that restriction of $c$ (resp. $c'$)
to $Q$ is equal to
the composite $a\circ c_Q$ (resp. $a'\circ c_Q$), and
$\theta\circ c=\theta'\circ c'$, with notation as in $(\clubsuit)$.
Then charts $c$ and $c'$ induce an isomorphism $\AAA \to \AAA'$
on $U$, which makes the diagram $(\clubsuit)$ commutative.

Next we shall prove the existence of a triple
$(\AAA,a,\theta)$. For a geometric point $\bar{s}$ on $S$,
consider the localization $S'=\Spec \OOO_{S,\bar{s}}$.
Set $Q=\alpha^{-1}f^{-1}\bar{\MMM}_{\Delta,\bar{s}}$.
Then by the case of the spectrum of a strictly Henselian local $k$-rings,
there exist a log structure $\AAA$ on $S'$,
a $\Delta^0$-FR morphism $a:\alpha^*f^*\MMM_{\Delta,\bar{s}}\to\AAA$, and the diagram of fine log structures on $S'$
\[
\xymatrix{
Q\ar[r]^{\bar{a}} \ar[d]^{c}& F \ar[r]^{\bar{\theta}} \ar[d]^{c'} & \bar{\NNN}_{\bar{s}}\ar[d]^{c''} \\
\alpha^*f^*\MMM_{\Delta,\bar{s}} \ar[r]^(0.7){a} & \AAA \ar[r]^{\theta} & \NNN_{\bar{s}}\\
}
\]
such that $\theta\circ a=h\circ h_f$.
Here $c$, $c'$, and $c''$ are charts and $F=\bar{\AAA}$.
To prove the existence on $S$, by the uniqueness,
it suffices only to show that we can extend the above diagram
to some \'etale neighborhood of $\bar{s}$.
In some \'etale neighborhood $U$ of $\bar{s}$,
there exists charts $\tilde{c}:Q\to \alpha^*f^*\MMM_{\Delta}$ and
$\tilde{c''}:\bar{\NNN}_{\bar{s}}\to \NNN$
extending $c$ and $c''$ respectively, such that
the diagram
\[
\xymatrix{
Q\ar[r]^{\bar{a}} \ar[d]^{\tilde{c}}& F \ar[r]^{\bar{\theta}}  & \bar{\NNN}_{\bar{s}}\ar[d]^{\tilde{c}''} \\
\alpha^*f^*\MMM_{\Delta} \ar[rr]^{h\circ h_f} &  & \NNN\\
}
\]
commutes.
Let $\tilde{\AAA}$ be the fine log structure
associated to the prelog structure $F\stackrel{\bar{\theta}}{\to} \bar{\NNN}_{\bar{s}}\stackrel{\tilde{c}''}{\to} \NNN \to \OOO_{U}$.
Then
there exists a sequence of morphisms of log structures
\[
\alpha^*f^*\MMM_{\Delta}\stackrel{\tilde{a}}{\longrightarrow} \tilde{\AAA} \stackrel{\tilde{\theta}}{\longrightarrow} \NNN,
\]
such that $\tilde{a}\circ\tilde{\theta}=h\circ h_f$,
which is an extension of
$\alpha^*f^*\MMM_{\Delta,\bar{s}}\stackrel{a}{\to} \AAA \stackrel{\theta}{\to} \NNN_{\bar{s}}$.
Since $\alpha^*f^*\MMM_{\Delta}\to \tilde{\AAA}$ has a chart by $Q\to F$,
thus by Lemma~\ref{sub1}
we conclude that $(f\circ \alpha,\tilde{a}:\alpha^*f^*\MMM_{\Delta}\to \tilde{\AAA}):(S,\tilde{\AAA}) \to (X_{\Delta},\MMM_{\Delta})$ is
a $\Delta^0$-FR morphism.
This completes the proof of Proposition~\ref{const}.
\QED

{\it Proof of Theorem 1.2.}
Let $\mathcal{H}om (\XXX_{(\Sigma,\Sigma^0)},\XXX_{\Delta,\Delta^0)})$
 be the category of torus-equivariant 1-morphisms
 from $\XXX_{(\Sigma,\Sigma^0)}$ to $\XXX_{(\Delta,\Delta^0)}$
 (whose morphisms are 2-isomorphisms).
 Let $\Hom ((\Sigma,\Sigma^0),(\Delta,\Delta^0))$ be
 the discrete category arising from
 the set of morphisms from $(\Sigma,\Sigma^0)$ to
 $(\Delta,\Delta^0)$.
We have to show that the natural map
\[
\Phi:\mathcal{H}om(\XXX_{(\Sigma,\Sigma^0)},\XXX_{(\Delta,\Delta^0)})
\to
\Hom ((\Sigma,\Sigma^0),(\Delta,\Delta^0)),
\]
is an equivalence.
This amounts to the following statement:
If $F:(\Sigma,\Sigma^0)\to (\Delta,\Delta^0)$ is a map of stacky fans
and 
$(f,h_f):(X_{\Sigma},\MMM_{\Sigma})\to (X_{\Delta},\MMM_{\Delta})$
denotes the torus-equivariant morphism
(with the natural morphism of the log structures)
of toric varieties
induced by $F:\Sigma\to \Delta$,
then
there exists a torus-equivariant
1-morphism (with the natural morphism of log structures)
\[
(\tilde{f},h_{\tilde{f}}):(\XXX_{(\Sigma,\Sigma^0)},\MMM_{(\Sigma,\Sigma^0)}) \to (\XXX_{(\Delta,\Delta^0)},\MMM_{(\Delta,\Delta^0)})
\]
such that
$(f,h_f)\circ(\pi_{(\Sigma,\Sigma^0)},h_{(\Sigma,\Sigma^0)})
=(\pi_{(\Delta,\Delta^0)},h_{(\Delta,\Delta^0)})\circ (\tilde{f},h_{\tilde{f}})$,
and it is unique up to a unique isomorphism.
By Proposition~\ref{const},
for each object $(\alpha,h):(S,\NNN)\to (X_{\Sigma},\MMM_{\Sigma})$
in $\XXX_{\Sigma}(\Sigma^0)$,
we can choose a pair $((f\circ \alpha,g),\xi_{(\alpha,h)})$
where $(f\circ \alpha,g):(S,\MMM)\to (X_{\Delta},\MMM_{\Delta})$
is an object
in $\XXX_{\Delta}(\Delta^0)$, i.e., a $\Delta^0$-FR morphism, and $\xi_{(\alpha,h)}:\MMM\to \NNN$ is a homomorphism of log structures
such that the diagram
\[
\xymatrix{
 \alpha^*f^*\MMM_{\Delta} \ar[r]^{g} \ar[d]_{\alpha^*h_f} & \MMM \ar[d]^{\xi_{(\alpha,h)}}\\
 \alpha^*\MMM_{\Sigma} \ar[r]^{h} & \NNN \\
 }
\]
commutes.
For each object $(\alpha,h):(S,\NNN)\to (X_{\Sigma},\MMM_{\Sigma})
\in \textup{Ob}(\XXX_{\Sigma}(\Sigma^0))$,
we choose such a pair $(\tilde{f}((\alpha,h)):(S,\MMM)\to (X_{\Delta},\MMM_{\Delta})\in \textup{Ob}(\XXX_{\Delta}(\Delta^0)),\xi_{(\alpha,h)}:\MMM\to \NNN)$.
By the axiom of choice, there exists a function
$\textup{Ob}(\XXX_{\Sigma}(\Sigma^0)) \to \textup{Ob}(\XXX_{\Delta}(\Delta^0))$, $(\alpha,h)\mapsto \tilde{f}((\alpha,h))$.
Let $(\alpha_i,h_i):(S_i,\NNN_i)\to (X_{\Sigma},\MMM_{\Sigma})$
be a $\Sigma^0$-FR morphism for $i=1,2$.
For each morphism $(q,e):(S_1,\NNN_1)\to (S_2,\NNN_2)$
in $\XXX_{\Sigma}(\Sigma^0)$,
define $\tilde{f}((q,e)):\tilde{f}((S_1,\NNN_1)):=(S_1,\MMM_1)_{/(X_{\Delta},\MMM_{\Delta})}\to
\tilde{f}((S_2,\NNN_2)):=(S_2,\MMM_2)_{/(X_{\Delta},\MMM_{\Delta})}$
to be $(q,\tilde{f}(e)):(S_1,\MMM_1)\to (S_2,\MMM_2)$
such that the diagram
\[
\xymatrix@R=4mm @C=10mm{
  & \MMM_1 \ar[rr]^{\xi_{(\alpha_1,h_1)}}&  & \NNN_1 \\
 q^*\MMM_2\ar[rr]^(0.7){q^*\xi_{(\alpha_2,h_2)}}\ar[ur]^{\tilde{f}(e)} &  & q^*\NNN_2\ar[ur]^{e} &  \\
  & \alpha_1^*f^*\MMM_{\Delta}\ar[uu]_(0.2){g_1} \ar[ul]^{q^*g_2}\ar[rr]&  & \alpha_1^*\MMM_{\Sigma}\ar[uu]^{h_1} \ar[lu]^{q^*h_2}\\
}
\]
commutes.
Here $\xi_{(\alpha_i,h_i)}$'s are the homomorphisms chosen as above.
The uniqueness of such a homomorphism $\tilde{f}(e)$
follows from Proposition~\ref{const}.
This yields a functor $\tilde{f}:\XXX_{\Sigma}(\Sigma^0)\to
\XXX_{\Delta}(\Delta^0)$ with a homomorphism
of log structures $\xi:\tilde{f}^*\MMM_{(\Delta,\Delta^0)}\to \MMM_{(\Sigma,\Sigma^0)}$ (it is determined by the collection $\{\xi_{(\alpha,h)}\}$).
It gives rise to a lifted morphism $(\tilde{f},\xi):
(\XXX_{(\Sigma,\Sigma^0)},\MMM_{(\Sigma,\Sigma^0)}) \to (\XXX_{(\Delta,\Delta^0)},\MMM_{(\Delta,\Delta^0)})$.
Since $\MMM_{(\Sigma,\Sigma^0)} \subset \OOO_{\XXX_{(\Sigma,\Sigma^0)}}$
on $\XXX_{(\Sigma,\Sigma^0),\textup{et}}$,
thus $\xi:\tilde{f}^*\MMM_{(\Delta,\Delta^0)}\to \MMM_{(\Sigma,\Sigma^0)}$
is the homomorphism $h_{\tilde{f}}$ induced by $\tilde{f}$ (cf. Section 2.4).
In addition,
if $\tilde{f}':\XXX_{(\Sigma,\Sigma^0)}\to \XXX_{(\Delta,\Delta^0)}$
is another lifting of $f$
and $\zeta:\tilde{f}\to \tilde{f}'$ is a 2-isomorphism,
then $\zeta$ induces an isomorphism
$\sigma:\tilde{f}^{*}\MMM_{(\Delta,\Delta^0)}
\to \tilde{f}^{'*}\MMM_{(\Delta,\Delta^0)}$
such that $h_{\tilde{f}}=\xi=h_{\tilde{f}'}\circ \sigma$,
where $h_{\tilde{f}'}:\tilde{f}^{'*}\MMM_{(\Delta,\Delta^0)}\to \MMM_{(\Sigma,\Sigma^0)}$ is the homomorphism induced by $\tilde{f}'$.
Therefore, for another lifting $\tilde{f}':\XXX_{(\Sigma,\Sigma^0)}=\XXX_{\Sigma}(\Sigma^0)\to \XXX_{\Delta}(\Delta^0)=\XXX_{(\Delta,\Delta^0)}$ of $f$,
the existence and the uniqueness of 2-isomorphism
$\zeta:\tilde{f}\to\tilde{f}'$ follows from Proposition~\ref{const}.
Indeed,
let $(\alpha,h):(S,\NNN)\to (X_{\Sigma},\MMM_{\Sigma})$
be a $\Sigma^0$-FR morphism and set
$\tilde{f}((\alpha,h))=\{ (f\circ \alpha, g):(S,\MMM)\to (X_{\Delta},\MMM_{\Delta})\}$
and 
$\tilde{f}'((\alpha,h))=\{ (f\circ \alpha, g'):(S,\MMM')\to (X_{\Delta},\MMM_{\Delta})\}$
(these are $\Delta^0$-FR morphisms).
Then we have the following commutative diagram
\[
\xymatrix@R=3mm @C=12mm{
 \alpha^*f^*\MMM_{\Delta} \ar[dd]^{\alpha^*h_f} \ar[rr]^{g} \ar[rd]^{g'} &   & \MMM \ar[dd]^{h_{\tilde{f}}}\\
  & \MMM' \ar[dr]^{h_{\tilde{f}'}}  &  \\
 \alpha^*\MMM_{\Sigma} \ar[rr]^{h} &   & \NNN \\
}
\]
where $h_{\tilde{f}}:\MMM\to \NNN$ (resp. $h_{\tilde{f}'}:\MMM'
\to \NNN$) denotes the homomorphism induced by $h_{\tilde{f}}:\tilde{f}^*\MMM_{(\Delta,\Delta^0)}\to \MMM_{(\Sigma,\Sigma^0)}$
(resp. $h_{\tilde{f}'}:\tilde{f}^{'*}\MMM_{(\Delta,\Delta^0)}\to \MMM_{(\Sigma,\Sigma^0)}$) (we abuse notation).
By Proposition~\ref{const},
there exists a unique isomorphism
of log structures $\sigma_{(\alpha,h)}:\MMM\to \MMM'$
which fits into the above diagram.
Then we can easily see that the collection
$\{ \sigma_{(\alpha,h)} \}_{(\alpha,h)\in\XXX_{\Sigma}(\Sigma^0)}$
defines a 2-isomorphism $\tilde{f}\to \tilde{f}'$.
Conversely, by the above observation and Proposition~\ref{const},
a 2-isomorphism $\tilde{f}\to \tilde{f}'$ must be
$\{ \sigma_{(\alpha,h)} \}_{(\alpha,h)\in\XXX_{\Sigma}(\Sigma^0)}$
and thus the uniqueness follows.
Finally, we shall show that $\tilde{f}$ is torus-equivariant
(cf. Definition~\ref{torusequiv}).
It follows from the uniqueness 
(up to a unique isomorphism) of
a lifting $\XXX_{(\Sigma,\Sigma^0)}\times \Spec k[M_1]\to \XXX_{(\Delta,\Delta^0)}$ of the torus-equivariant morphism $X_{\Sigma}\times\Spec k[M_1]\stackrel{a}{\to} X_{\Sigma}\stackrel{f}{\to}
X_{\Delta}$. Here  $a$ is the torus action and we mean by a lifting of $f\circ a$ a functor which commutes with $f\circ a$ via coarse moduli maps.
Thus we complete the proof of Theorem~\ref{main}.
\QED

Theorem~\ref{main} and its proof imply the followings:

\begin{Corollary}
\label{uniq}
Let $f:X_{\Sigma}\to X_{\Delta}$ be a torus-equivariant morphism
of simplicial toric varieties.
Then a functor $($not necessarily torus-equivariant$)$
$\tilde{f}:\XXX_{(\Sigma,\Sigma^0)}\to \XXX_{(\Delta,\Delta^0)}$
such that
$\pi_{(\Delta,\Delta^0)}\circ \tilde{f}=f\circ \pi_{(\Sigma,\Sigma^0)}$
is unique up to a unique isomorphism $($if it exists$)$.
\end{Corollary}
\Proof
It follows immediately from the proof of Theorem~\ref{main}.
\QED

\begin{Corollary}
\label{equivcha}
Let $f:\XXX_{(\Sigma,\Sigma^0)}\to \XXX_{(\Delta,\Delta^0)}$ be
a $($not necessarily torus-equivariant$)$ functor.
Then $f$ is torus-equivariant if and only if
the induced morphism
$f_c:X_{\Sigma}\to X_{\Delta}$ of toric varieties is torus-equivariant.

\end{Corollary}

\Proof
The ``only if" part follows from Proposition~\ref{eqeq}.
The proof of Theorem~\ref{main} implies the ``if" part.
\QED

\begin{Corollary}
Let $\Sigma$ and $\Delta$ be simplicial fans and
let $(\Delta,\Delta^0)$ be a stacky fan that is an extension of $\Delta$.
Let $F:\Sigma\to \Delta$  be a homomorphism of fans and
let $f:X_{\Sigma}\to X_{\Delta}$ be the associated morphism
of toric varieties. Then there exist a stacky fan
$(\Sigma,\Sigma^0)$ that is an extension of $\Sigma$
and a torus-equivariant morphism $\tilde{f}:\XXX_{(\Sigma,\Sigma^0)}\to
\XXX_{(\Delta,\Delta^0)}$ such that $\pi_{(\Delta,\Delta^0)}\circ \tilde{f}=
f\circ \pi_{(\Sigma,\Sigma^0)}$.
Moreover if we fix such a stacky fan $(\Sigma,\Sigma^0)$,
then $\tilde{f}$ is unique up to a unique isomorphism.
\end{Corollary}
\Proof
Theorem~\ref{main} immediately implies our assertion
because we can choose a free-net $\Sigma^0$ such that $F(\Sigma^0)\subset \Delta^0$.
\QED

\begin{Corollary}
\label{cor}
Let $(\Sigma,\Sigma^0)$ be a stacky fan in $N_{\RR}$
and let $\XXX_{(\Sigma,\Sigma^0)}$ be the associated toric algebraic stack.
Then there exists a smooth surjective torus-equivariant morphism
\[
p:X_{\Delta} \longrightarrow \XXX_{(\Sigma,\Sigma^0)}
\]
where $X_{\Delta}$ is a quasi-affine smooth toric variety.
Furthermore, $X_{\Delta}$ can be explicitly constructed.
\end{Corollary}

\Proof
Without loss of generality, we may suppose that rays of $\Sigma$ span the vector space $N_{\RR}$.
Set $\tilde{N}=\oplus_{\rho\in \Sigma(1)}\ZZ\cdot e_{\rho}$.
Define a homomorphism of abelian groups
$\eta:\tilde{N}\longrightarrow N$
by $e_{\rho}\mapsto P_{\rho}$ where $P_{\rho}$ is the generator of
$\Sigma^0$ on $\rho$ (cf. Definition~\ref{stfan}).
Let $\Delta$ be a fan in $\tilde{N}_{\RR}$
that consists of cones  $\gamma$ such that
$\gamma$ is a face of the cone $\oplus_{\rho\in \Sigma(1)}\RR_{\ge0}\cdot e_{\rho}$ and
$\eta_{\RR}(\gamma)$ lies in $\Delta$.
If $\Delta_{\extup}^0$ denotes the canonical free-net
(cf. Definition~\ref{stfan}), then $\eta$ induces
the morphism of stacky fans $\eta:(\Delta,\Delta_{\extup}^0)
\to (\Sigma,\Sigma^0)$.
Note that $\XXX_{(\Delta,\Delta_{\extup}^0)}$
is the quasi-affine smooth toric variety $X_{\Delta}$.
Let $p:X_{\Delta}\to \XXX_{(\Sigma,\Sigma^0)}$ be
the torus-equivariant morphism induced by $\eta$ (cf. Theorem~\ref{main}).
Since the composite $q:=\pi_{(\Sigma,\Sigma^0)}\circ p:X_{\Delta}\to\XXX_{(\Sigma,\Sigma^0)}\to X_{\Sigma}$ is surjective and $\pi_{(\Sigma,\Sigma^0)}$ is
the coarse moduli map, thus by \cite[Proposition 5.4 (ii)]{LM},
$p$ is surjective. It remains to show that $p$ is smooth.
This is an application of K. Kato's notion of {\it log smoothness}.
From the construction of $\eta$, Lemma~\ref{sub1} and Lemma~\ref{constfreeres}, we can easily see that
the induced morphism
$(q,h_q):(X_{\Delta},\MMM_{\Delta})\to(X_{\Sigma},\MMM_{\Sigma})$
is a $\Sigma^0$-FR morphism.
Moreover, by \cite[Theorem 3.5]{log}, we see that $(q,h_q)$ is log smooth
($\textup{ch}(k)=0$).
By the modular interpretation of $\XXX_{\Sigma}(\Sigma^0)$ (cf. Section 2),
there exists the 1-morphism $p':X_{\Delta}\to \XXX_{(\Sigma,\Sigma^0)}$
which corresponds to $(q,h_q)$. Theorem~\ref{main} and Corollary~\ref{uniq}
imply that
$p'$ coincides with $p$, and thus the
morphism $(p,h_p):(X_{\Delta},\MMM_{\Delta})\to (\XXX_{(\Sigma,\Sigma^0)},\MMM_{(\Sigma,\Sigma^0)})$
is a strict morphism.
Here $h_p$ is the homomorphism induced by $p$.
Then the following lemma implies that $p$ is smooth.

\begin{Lemma}
Let $(X,\MMM)$ be a log scheme and $(f,h):(X,\MMM)\to (\XXX_{(\Sigma,\Sigma^0)},\MMM_{(\Sigma,\Sigma^0)})$
a strict 1-morphism. If the composite
\[
(\pi_{(\Sigma,\Sigma^0)},h_{(\Sigma,\Sigma^0)})\circ (f,h):(X,\MMM)\to(X_{\Sigma},\MMM_{\Sigma})
\]
is formally log smooth, then $f:X\to \XXX_{(\Sigma,\Sigma^0)}$ is formally smooth.

\end{Lemma}

\Proof
It suffices to show the lifting property as in \cite[Definition 4.5]{OL}.
Let $i:T_0\to T$ be a closed immersion of schemes defined by a
square zero ideal.
Let $a_0:T_0\to X$ and $b:T\to \XXX_{(\Sigma,\Sigma^0)}$
be a pair of 1-morphisms such that $f\circ a_0\cong b\circ i$.
We have to show that there exists a 1-morphism
$a:T\to X$ such that $a\circ i=a_0$ and $f\circ a\cong b$.
There exists the following commutative diagram
\[
\xymatrix{
(T_0,a_0^*\MMM)\ar[d]^{i}\ar[r]^{a_0} &(X,\MMM) \ar[d]^{(f,h)}\\
(T,b^*\MMM_{(\Sigma,\Sigma^0)}) \ar[r]^(0.4){b} \ar[rd]& (\XXX_{(\Sigma,\Sigma^0)},\MMM_{(\Sigma,\Sigma^0)}) \ar[d]^{(\pi_{(\Sigma,\Sigma^0)},h_{(\Sigma,\Sigma^0)})}\\
  &  (X_{\Sigma},\MMM_{\Sigma}), \\
}
\]
where $i$, $b$, and $a_0$ denote induced
strict morphisms (we abuse notation).
Since $(\pi_{(\Sigma,\Sigma^0)},h_{(\Sigma,\Sigma^0)})\circ (f,h)$
is formally log smooth, there exists a morphism
\[
(a,v):(T,b^*\MMM_{(\Sigma,\Sigma^0)})\to (X,\MMM)
\]
such that $a_0=a\circ i$.
Thus it suffices only to prove $b\cong f\circ a$.
This is equivalent to showing that
$(a,v)$ is a strict morphism
since $\XXX_{\Sigma}(\Sigma^0)=\XXX_{(\Sigma,\Sigma^0)}$ is the moduli stack
of $\Delta^0$-FR morphisms into $(X_{\Sigma},\MMM_{\Sigma})$.
To see this, 
we have to show that for any geometric point $\bar{t}\to T$,
$\bar{v}_{\bar{t}}:(a^{-1}\bar{\MMM})_{\bar{t}}\to (b^*\MMM_{(\Sigma,\Sigma^0)}/\OOO_{T}^*)_{\bar{t}}$
is an isomorphism.
It follows from the following:
Let $\iota:P\to \NN^r$ be an injective homomorphism
of monoids such that $\iota(P)$ is close to $\NN^r$.
Let $e:\NN^r\to \NN^r$ be an endomorphism such that $\iota=e\circ \iota$.
Then $e$ is an isomorphism.
\QED

\section{A geometric characterization theorem}

The aim of this Section is to give proofs of Theorem~\ref{main2}
and Theorem~\ref{main3}. In this Section, we work over
an algebraically closed base field $k$ of characteristic zero,
except in Lemma~\ref{normal}.

\begin{Lemma}
\label{normal}
Let $\SSS$ be a normal Deligne-Mumford stack locally of finite type and
separated over a locally noetherian scheme,
and let $p:\SSS\to S$ be a coarse moduli map.
Then $S$ is normal.
\end{Lemma}

\Proof
Our assertion is \'etale local on $S$, and thus
we may assume that $S$ is the spectrum of a strictly Henselian local ring.
Set $S=\Spec O$.
In this situation, by \cite[Lemma 2.2.3]{AV}
there exist a normal strictly Henselian
local ring $R$, a finite group $G$ and an
action $m:\Spec R\times G\to \Spec R$
such that the quotient stack $[\Spec R/G]$
is isomorphic to $\SSS$,
and $R^{G}=O$ (here $R^G$ is the invariant ring).
Let $n:\Spec A\to \Spec O$ be the normalization of $\Spec O$
in the function field $\QQ(O)$.
Let us denote by $q:\Spec R\to \Spec O$
(resp. $\textup{pr}_1:\Spec R\times G\to \Spec R$)
the composite $\Spec R\to \SSS\to S=\Spec O$
(resp. the natural projection).
Then by the universality of the normalization,
there exists a unique morphism $\tilde{q}:\Spec R\to \Spec A$
such that $n\circ \tilde{q}=q$.
Note that
$\tilde{q}\circ \textup{pr}_1$  (resp. $\tilde{q}\circ m$)
is the unique lifting of
$q\circ \textup{pr}_1$ (resp. $q\circ m$).
Since $q\circ \textup{pr}_1=q\circ m$, thus
we have $\tilde{q}\circ \textup{pr}_1=\tilde{q}\circ m$.
This implies that $A\subset O=R^G$ and we conclude that
$S$ is normal.
\QED

\begin{Proposition}
\label{setup}
Let $(\mathcal{X}, \iota:\GG_m^d\hookrightarrow \mathcal{X}, a:\mathcal{X}\times\GG_m^d\to \mathcal{X})$ be a toric triple over $k$.
Then the complement
$\DDD:=\mathcal{X}-\GG_m^d$ with reduced closed substack structure
is a divisor with normal crossings, and
the coarse moduli space $X$ is a simplicial toric variety over $k$.
\end{Proposition}

\Proof
First, we shall prove that $X$ is a toric variety over $k$.
Observe that the coarse moduli scheme
$X$ is a normal variety over $k$, i.e.,
normal and of finite type and separated over $k$.
Indeed, according to Keel-Mori Theorem
$X$ is locally of finite type and separated over $k$.
Since $\XX$ is of finite type over $k$ and the underlying continuous
morphism $|\XX|\to |X|$ (cf. \cite[5]{LM})
of the coarse moduli map
 is a homeomorphism,
thus $X$ is of finite type over $k$ by \cite[5.6.3]{LM}.
Since $\XX$ is smooth over $k$, thus by Lemma~\ref{normal}
$X$ is normal.
Since  $\iota:\GG_m^d\hookrightarrow \mathcal{X}$, the coarse moduli space $X$ contains
$\GG_m^d$ as a dense open subset. 
The torus action $\mathcal{X}\times\GG_m^d \to \mathcal{X}$,
gives rise to a morphism of coarse moduli spaces
$a_0:X\times \GG_m^d\to X$ because $X\times \GG_m^d$ is
a coarse moduli space for $\mathcal{X}\times\GG_m^d$.
Moreover by the universality of coarse moduli spaces,
it is an action of $\GG_m^d$ on $X$.
Therefore $X$ is a toric variety over $k$.

Next we shall prove that
the complement
$\DDD$ is a divisor with normal crossings.
Set $\GG_m^d=\Spec k[M]$ ($M=\ZZ^d$) and $X_{\Sigma}=X$
where $\Sigma$
is a fan in $N\otimes_{\ZZ}\RR$
($N=\Hom_{\ZZ}(M,\ZZ))$.
Let $\bar{x}\to X_{\Sigma}$ be a geometric point on $X_{\Sigma}$ and put $\OOO:=\OOO_{X_{\Sigma},\bar{x}}$ (the \'etale stalk).
Consider the pull-back $\mathcal{X}_{\OOO}:=\mathcal{X}\times_{X_{\Sigma}}\OOO\to \Spec \OOO$ by $\Spec \OOO_{X_{\Sigma},\bar{x}}\to X_{\Sigma}$.
Clearly, 
our assertion is an \'etale local issue on $X_{\Sigma}$
and thus it suffices to show that $\DDD$ defines a divisor with normal
crossings on $\mathcal{X}_{\OOO}$.
By \cite[Lemma 2.2.3]{AV}, there exists a strictly Henselian local
$k$-ring $R$ and a finite group $\Gamma$ acting on $\Spec R$
such that $\mathcal{X}_{\OOO}\cong [\Spec R/\Gamma]$.
We have a sequence of morphisms
\[
\Spec R\stackrel{p}{\to} [\Spec R/\Gamma] \stackrel{\pi}{\to} \Spec \OOO.
\]
The composite $q:=\pi \circ p$ is a finite surjective morphism.
If $U$ denotes the open subscheme of $\Spec \OOO$ which is
induced by the torus embedding
$\GG_m^d\subset X_{\Sigma}$, then the restriction
$q^{-1}(U)\to U$
is a finite \'etale surjective morphism.
Let us denote by $\MMM_{\OOO}$ the pull-back of the canonical log
structure $\MMM_{\Sigma}$ on $X_{\Sigma}$ to $\Spec \OOO$.
Then in virtue of log Nagata-Zariski purity Theorem \cite[Theorem 3.3]{Mot} (See also \cite[Remark 1.10]{Ho}),
the complement $\Spec R-q^{-1}(U)$ (or equivalently $\DDD$)
defines a log structure on $\Spec R$ 
(we shall denote by $\MMM_R$ this log structure) and 
the finite \'etale surjective morphism $q^{-1}(U)\to U$ extends to
a Kummer log \'etale surjective morphism
\[
(q,h):(\Spec R, \MMM_R) \to (\Spec \OOO, \MMM_{\OOO}).
\]
Let $\hat{\OOO}$ be the completion of $\OOO$ along its maximal ideal.
Let us denote by
\[
(\hat{q},\hat{h}):(T,\MMM_{R}|_T) \to (\Spec \hat{\OOO}, \MMM_{\hat{\OOO}}:=\MMM_{\OOO}|_{\hat{\OOO}})
\]
the pull-back of $(q,h)$ by $\Spec \hat{\OOO}\to \Spec
\OOO$.
Then by \cite[Theorem 3.2]{log2}, the log scheme
$(\Spec \hat{\OOO}, \MMM_{\OOO}|_{\hat{\OOO}})$ is isomorphic to
\[
\Spec (P\to k(\bar{x})[[P]][[\NN^l]])
\]
where
$k(\bar{x})$ is the residue field of $\bar{x}\to X_{\Sigma}$,
$P:=\bar{\MMM}_{\hat{\OOO},\bar{x}}\to k(\bar{x})[[P]][[\NN^l]]\ (p\mapsto p)$,
and $l$ is a non-negative integer.
(Strictly speaking, \cite{log2} only treats the case of Zariski log structures,
but the same proof can apply to the case of \'etale log structures.)
By taking a connected component of $T$ if necessary,
we may assume that $T$ is connected.
Note that since the connected scheme $T$ is finite over $\Spec \hat{\OOO}$,
$T$ is the spectrum of a strictly Henselian local $k$-ring.
Let $Q$ be the stalk of $\bar{\MMM}_R|_T$ at a geometric
point $\bar{t}\to T$ lying over the closed point $t$ of $T$.
Then by \cite[Proposition A.4]{Ho},
the Kummer log \'etale cover $(\hat{q},\hat{h})$
has the form
\[
\xymatrix{
\Spec (Q\to \ZZ[Q]\otimes_{\ZZ[P]}k(\bar{x})[[P]][[\NN^l]]) \ar[r]
& \Spec (P\to k(\bar{x})[[P]][[\NN^l]]) \\
}
\]
defined by $P=\bar{\MMM}_{\hat{\OOO},\bar{x}}\to (\bar{\MMM}_{R}|_T)_{\bar{t}}=Q$,
$\textup{Id}_{\NN^l}:\NN^l\to \NN^l$
and the natural map
$Q\to \ZZ[Q]\otimes_{\ZZ[P]}k(\bar{x})[[P]][[\NN^l]]$.
(Note that $\ZZ[Q]\otimes_{\ZZ[P]}k(\bar{x})[[P]]\cong k(\bar{x})[[Q]]$
because $Q\to P$ is Kummer.)
Since $T$ is regular, thus $Q$ is free, i.e., $Q\cong \NN^{r}$ for some
non-negative integer $r$.
This implies that $\DDD$ is a divisor with normal crossings.

Finally, we shall show that the toric variety
$X_{\Sigma}$ is simplicial.
To this end, we assume that $\Sigma$ is not simplicial
and show that such an assumption gives rise to a contradiction.
From the assumption, there exists a geometric point $\alpha:\bar{x}\to X_{\Sigma}$
such that the number of irreducible components of the complement $D:=X_{\Sigma}-\GG_m^d$ on which
the point $\bar{x}$ lies is greater than
the rank $\rank\bar{\MMM}_{\Sigma,\bar{x}}^{gp}$
of $\bar{\MMM}_{\Sigma,\bar{x}}^{gp}$.
Let $r$ be the number of irreducible components
of $D$ on which the point $\bar{x}$ lies.
Put $P:=\bar{\MMM}_{\Sigma,\bar{x}}$.
By the same argument as above,
there exist a strictly Henselian local $k$-ring $R$ and a sequence
of Kummer log \'etale covers
\[
(\Spec R, p^*\MMM_{\DDD})\stackrel{p}{\to} (\mathcal{X}\times_{X_{\Sigma}}\Spec\OOO_{X_{\Sigma},\bar{x}},\MMM_{\DDD})\to (\Spec \OOO_{X_{\Sigma},\bar{x}},\MMM_{\Sigma}|_{\Spec \OOO_{X_{\Sigma},\bar{x}}})
\]
where $\MMM_{\DDD}$ is the log structure induced by $\DDD$
and the left morphism is a strict morphism.
Moreover the pull-back of the composite
$\Spec R\to \Spec\OOO_{\Sigma,\bar{x}}$ (it is a finite morphism)
by the completion
$\Spec \hat{\OOO}_{\Sigma,\bar{x}}\to \Spec\OOO_{\Sigma,\bar{x}}$
along the maximal ideal is of the form
\[
\Spec k(\bar{x})[[\NN^r]][[\NN^l]] \to \Spec k(\bar{x})[[P]][[\NN^l]]
\]
because $\pi^{-1}(D)_{\xtup}=\DDD$
and $\DDD$ is a normal crossing divisor on the smooth stack.
Here $\pi:\mathcal{X}\to X_{\Sigma}$ is the coarse moduli map,
and for each irreducible component $C$ of $D$,
$\pi^{-1}(C)_{\xtup}$ is an irreducible component of $\DDD$
because the underlying continuous map
$|\pi|:|\XX|\to |X|$ (cf. \cite[5.2]{LM}) is a homeomorphism.
However we have 
\[
\dim \Spec k(\bar{x})[[\NN^r]][[\NN^l]]>
\dim\Spec k(\bar{x})[[P]][[\NN^l]]=\rank P^{gp}+l.
\]
This is a contradiction.
\QED

{\it Proof of Theorem~\ref{main2}.}
We shall construct a morphism
from $\XX$ to some toric algebraic stack and show that it is an isomorphism
with desired properties.

\noindent
(Step 1)
We first construct a morphism from $\XX$ to some toric algebraic stack.
Set $\GG_m^d=\Spec k[M]$ ($M=\ZZ^d$) and $N=\Hom_{\ZZ}(M,\ZZ)$.
By Proposition~\ref{setup},
we can put $X=X_{\Sigma}$ where $\Sigma$ is a simplicial fan
in $N\otimes_{\ZZ}\RR$,
and let us denote by $\pi:\XX\to X_{\Sigma}$ the coarse moduli map.
By \cite[3.3 (2)]{I2} and \cite[4.4]{I2},
there exists a morphism $\phi:\XX\to \XXX_{(\Sigma,\Sigma_{\extup}^0)}$
such that $\pi\cong \pi_{(\Sigma,\Sigma_{\extup}^0)}\circ \phi$.
For a ray $\rho\in \Sigma(1)$, we denote by $V(\rho)$
the corresponding irreducible
component of $D=X_{\Sigma}-\Spec k[M]$
where $\Spec k[M]\subset X_{\Sigma}$ is the torus embedding,
that is, the torus-invariant divisor corresponding to $\rho$.
Then $\VVV(\rho):=\pi_{(\Sigma,\Sigma_{\extup}^0)}^{-1}(V(\rho))_{\xtup}$ (resp. $\WWW(\rho):=\pi^{-1}(V(\rho))_{\xtup}$) is
an irreducible component of the normal crossing divisor
$\XXX_{(\Sigma,\Sigma_{\extup}^0)}-\Spec k[M]$
(resp. $\DDD=\XX-\GG_m^d$).
Since $\pi_{(\Sigma,\Sigma_{\extup}^0)}$ and $\pi$
are coarse moduli maps,
$\phi^{-1}(\VVV(\rho))_{\xtup}=\WWW(\rho)$.
For each ray $\rho\in \Sigma(1)$, let $n_{\rho}\in \NN$ be the natural
number such that
\[
\phi^{-1}(\VVV(\rho))=n_{\rho}\cdot \WWW(\rho).
\]
Let $(\Sigma,\Sigma^0)$ be the stacky fan whose level on
each ray $\rho$
is $n_{\rho}$.
If $\MMM_{\DDD}$ denotes the log structure associated to $\DDD$,
the morphism of log stacks $(\pi,h_{\pi}):(\XX,\MMM_{\DDD})\to (X_{\Sigma},\MMM_{\Sigma})$ (cf. Section 2.4) is a $\Sigma^0$-FR morphism
since $\DDD$ is a normal crossing divisor and $(\pi_{(\Sigma,\Sigma_{\extup}^0)},h_{(\Sigma,\Sigma_{\extup}^0)}):(\XXX_{(\Sigma,\Sigma_{\extup}^0)},\MMM_{(\Sigma,\Sigma_{\extup}^0)})\to (X_{\Sigma},\MMM_{\Sigma})$
is a $\Sigma_{\extup}^0$-FR morphism.
Then there exists a strict morphism of log stacks
\[
\Phi:(\XX,\MMM_{\DDD})\longrightarrow (\XXX_{(\Sigma,\Sigma^0)},\MMM_{(\Sigma,\Sigma^0)})
\]
over $(X_{\Sigma},\MMM_{\Sigma})$,
which is associated to the $\Sigma^0$-FR morphism $(\pi,h_{\pi})$.
By the construction of $\Phi$,
the restriction of $\Phi$
to $\GG_m^d\subset \XX$
induces an isomorphism $\GG_m^d\stackrel{\sim}{\to} \Spec k[M]\subset \XXX_{(\Sigma,\Sigma^0)}$
of group $k$-schemes.
\vspace{1mm}

We will prove that $\Phi$ is an isomorphism in (Step 2) and (Step 3).

\vspace{1mm}

\noindent
(Step 2) Observe that it suffices to prove that for
each closed point $\bar{x}:=\Spec k\to X_{\Sigma}$,
the pull-back $\Phi_{\bar{x}}:\XX\times_{X_{\Sigma}}\Spec \OOO_{X_{\Sigma},\bar{x}} \to \XXX_{(\Sigma,\Sigma^0)}\times_{X_{\Sigma}}\Spec \OOO_{X_{\Sigma},\bar{x}}$ by $\Spec \OOO_{X_{\Sigma},\bar{x}}\to X_{\Sigma}$ is an isomorphism.
(Here $\OOO_{X_{\Sigma},\bar{x}}$ is the \'etale stalk.)
Indeed, assume that $\Phi_{\bar{x}}$ is an isomorphism for every closed point
$\bar{x}=\Spec k\to X_{\Sigma}$.
Then by \cite[Theorem 2.2.5]{B}, $\Phi$ is representable.
Moreover $\Phi$ is finite.
We see this as follows: Note that $\XX$ is separated over $k$,
thus $\Phi$ is separated.
In addition, clearly, $\Phi$ is of finite type.
Since $\pi$ and $\pi_{(\Sigma,\Sigma^0)}$ are coarse moduli maps (in particular, proper),
thus by \cite[Proposition 2.7]{hom},
$\Phi$ is a proper and quasi-finite surjective morphism, i.e.,
a finite
surjective morphism 
(cf. \cite[Corollary A.2.1]{LM}).
It is an \'etale local issue on $X_{\Sigma}$
whether or not $\Phi$ is an isomorphism, and thus by
\cite[Chap. IV 8.8.2.4]{EGA} we conclude that $\Phi$ is an isomorphism
because $\Phi$ is a finite representable morphism.
Therefore, we shall prove that
$\Phi_{\bar{x}}:\XX\times_{X_{\Sigma}}\Spec \OOO_{X_{\Sigma},\bar{x}} \to \XXX_{(\Sigma,\Sigma^0)}\times_{X_{\Sigma}}\Spec \OOO_{X_{\Sigma},\bar{x}}$
is an isomorphism for each closed point $\bar{x}=\Spec k\to X_{\Sigma}$.
For simplicity, put
$\OOO:=\OOO_{X_{\Sigma},\bar{x}}$,
$\XX':=\XX\times_{X_{\Sigma}}\Spec \OOO_{X_{\Sigma},\bar{x}}$
and
$\XXX'_{(\Sigma,\Sigma^0)}:=\XXX_{(\Sigma,\Sigma^0)}\times_{X_{\Sigma}}\Spec \OOO_{X_{\Sigma},\bar{x}}$.
Set $\alpha:\Spec \OOO\to X_{\Sigma}$.
Write $\MMM$,
$\MMM_{\DDD}'$, and $\MMM_{(\Sigma,\Sigma^0)}'$
for log structures $\alpha^*\MMM_{\Sigma}$,
$(\alpha\times_{X_{\Sigma}}\XX)^*\MMM_{\DDD}$
and $(\alpha\times_{X_{\Sigma}}\XXX_{(\Sigma,\Sigma^0)})^*\MMM_{(\Sigma,\Sigma^0)}$ on $\Spec \OOO$, $\XX'$ and $\XXX_{(\Sigma,\Sigma^0)}'$ respectively.
Clearly, we may assume that $\Sigma$ is a simplicial cone
$\sigma$ and $X_{\Sigma}=\Spec k[P]\times_{k}\GG_m^l$
where $P=(\sigma^{\vee}\cap M)/(\textup{invertible\ elements})$
and $l$ is a non-negative integer.
In addition, by replacing $\sigma$ with a face if
necessary, we can suppose that the closed point $\bar{x}$ lies
on the torus orbit of the point $(o,1)\in \Spec k[P]\times_{k}\GG_m^l$.
Here $o\in \Spec k[P]$ is the origin, and $1\in \GG_m^l$
is the unit point.
Thus we may assume that $\bar{x}=(o,1)$.

\noindent
(Step 3)
We will prove that $\Phi_{\bar{x}}$ is an isomorphism.
To this end,
we first give an explicit representation of $(\XX',\MMM_{\DDD}')$
as a form of quotient stack.
By \cite[Lemma 2.2.3]{AV} and \cite[Theorem 2.12]{hom},
there exist a $d$-dimensional strictly Henselian regular local $k$-ring $R$
(here $d:=\dim X_{\Sigma}$),
a finite group $\Gamma$ acting on $R$ which is isomorphic
to the stabilizer group of any geometric point on $\XX$ lying
over $\bar{x}$,
and an isomorphism 
\[
\XX'\cong [\Spec R/\Gamma]
\]
over $\Spec \OOO$.
Furthermore the action of $\Gamma$ on the closed point of $\Spec R$ is trivial
and the invariant ring $R^{\Gamma}$ is the image of $\OOO\hookrightarrow R$.
Note that
if $\Aut_{\scriptsize{\Spec} \OOO}(\Spec R)$ 
denotes the group of automorphisms of $\Spec R$ over $\Spec \OOO$,
the natural homomorphism of groups
$\Gamma \to \Aut_{\scriptsize{\Spec} \OOO}(\Spec R)$ 
is injective because $\XX'$ is generically representable.
Let us denote by $p:\Spec R\to [\Spec R/\Gamma]$
the natural projection and put $\MMM_{\DDD,R}':=p^*\MMM_{\DDD}'$.
Consider the composite $\Spec R\to [\Spec R/\Gamma]\to \Spec \OOO$.
Then this composite induces the morphism of log schemes
\[
(f,h):(\Spec R,\MMM_{\DDD,R}')\to (\Spec \OOO,\MMM)
\]
whose underlying morphism $\Spec R\to \Spec \OOO$ is
finite and surjective.
Let $W$ be the open subscheme $\alpha^{-1}(\GG_m^d)\subset \Spec \OOO$
$(\GG_m^d\subset X_{\Sigma}$).
Then the restriction $f^{-1}(W)\to W$ is a finite \'etale surjective morphism.
In virtue of log Nagata-Zariski purity Theorem (\cite[Theorem 3.3]{Mot}
and \cite[Remark 1.10]{Ho}), $(f,h)$ is a Kummer log \'etale cover
($\textup{ch}(k)=0$).
Now put $\OOO=k\{ P,t_1,\ldots,t_l\}\subset k[[P]][[t_1,\ldots,t_l]]$
where $k\{ P,t_1,\ldots,t_l\}$ is the (strict) Henselization of the
Zariski stalk of the origin of $\Spec k[P,t_1,\ldots,t_l]$.
Consider the homomorphism $P=\bar{\MMM}_{\bar{s}}\to F:=\bar{\MMM}_{\DDD,R,\bar{t}}'$
where $\bar{s}\to \Spec \OOO$ and $\bar{t}\to \Spec R$ are
geometric points lying
over the closed points of $\Spec \OOO$ and $\Spec R$
respectively.
Then by \cite[Proposition A.4]{Ho}, $(f,h)$ is of the form
\[
\xymatrix{
\Spec (F\to k[F]\otimes_{k[P]}k\{ P,t_1,\ldots,t_l\})  \ar[r] & \Spec (P\to k\{ P,t_1,\ldots,t_l\}) \\
}
\]
where underlying morphism and homomorphism
of log structures are naturally induced by $P\to F$ and $t_i\mapsto t_i$.
Here $F\to k[F]\otimes_{k[P]}k\{ P,t_1,\ldots,t_l\}$
and
$P\to k\{ P,t_1,\ldots,t_l\}$ are the natural homomorphisms.
As observed in \cite[3.19]{Stix},
the group $\Aut_{(\scriptsize{\Spec} \OOO,\MMM)}((\Spec R,\MMM_{\DDD,R}'))$ of automorphisms of
$(\Spec R,\MMM_{\DDD,R}')$ over $(\Spec \OOO,\MMM)$
is naturally isomorphic to
$G:=\Hom_{\scriptsize{\textup{group}}}(F^{gp}/P^{gp},k^*)$.
Here an element $g\in G$ acts on $k[F]\otimes_{k[P]}k\{ P,t_1,\ldots,t_l\}$
by $f\mapsto g(f)\cdot f$ for any $f\in F$.
The natural forgetting homomorphism $\Aut_{(\scriptsize{\Spec} \OOO,\MMM)}((\Spec R,\MMM_{\DDD,R}'))\to \Aut_{\scriptsize{\Spec} \OOO}(\Spec R)$ is an isomorphism.
The injectivity is clear from the action of $G$, and the surjectivity 
follows from the facts $\MMM_{\DDD,R}'=\{ s \in \OOO_{\scriptsize{\Spec} R}|\ s\ \textup{is invertible on}\ f^{-1}(W) \}$
and $\MMM=\{s\in \OOO_{\scriptsize{\Spec} \OOO}|\ s\ \textup{is invertible on}\ W \}$.
Furthermore, since the category of Kummer log
\'etale coverings is a Galois category (cf. \cite[Theorem A.1]{Ho}),
the injective morphism
$\Gamma\to G=\Aut_{\scriptsize{\Spec}\OOO}(\Spec R)=\Aut_{(\scriptsize{\Spec} \OOO,\MMM)}((\Spec R,\MMM_{\DDD,R}'))$
is surjective, i.e., bijective.
Indeed, if it is not surjective,
then the Kummer log \'etale cover (or its underlying morphism)
$(\Spec R,\MMM_{\DDD,R}')/\Gamma\to (\Spec \OOO,\MMM)$
is not an isomorphism and thus
we obtain a contradiction to $R^{\Gamma}=\OOO$.
Since any local ring that is finite over a Henselian local ring
is also Henselian, thus we have $k[F]\otimes_{k[P]}k\{ P,t_1,\ldots,t_l\}\cong
k\{ F,t_1\ldots,t_l\}$ where
$k\{ F,t_1,\ldots,t_l\}$ is the (strict) Henselization of the
Zariski stalk of the origin of $\Spec k[F,t_1,\ldots,t_l]$.
Hence
there exists an isomorphism of log stacks
\[
(\XX',\MMM_{\DDD}')\cong ([\Spec k\{ F,t_1,\ldots,t_l\}/G],
\MMM_F)
\]
over $(\Spec \OOO,\MMM)$
where $\MMM_F$ is
the log structure on $[\Spec k\{ F,t_1,\ldots,t_l\}/G]$
induced by the natural chart $F\to k\{F,t_1,\ldots,t_l\}$.
In particular, the morphism
\[
([\Spec k\{ F,t_1,\ldots,t_l\}/G], \MMM_F)\to (\Spec \OOO, \MMM)
\]
is isomorphic to $(\XX',\MMM_{\DDD}')$
as a $\Sigma^0$-FR morphism over $(\Spec \OOO,\MMM)$.
Next, by using this form we will prove that $\Phi_{\bar{x}}$
is an isomorphism.
Note that the morphism $\Phi_{\bar{x}}:\XX'\to \XXX_{(\Sigma,\Sigma^0)}'$
over $\Spec \OOO$ is
the morphism associated to the
$\Sigma^0$-FR morphism
$(\XX',\MMM_{\DDD}')\to (\Spec \OOO, \MMM)$.
Thus what we have to show
is that $[\Spec k\{ F,t_1,\ldots,t_l\}/G]_{/\scriptsize{\Spec} \OOO}$ is the stack whose objects over $S\to \Spec \OOO$
are 
$\Sigma^0$-FR morphisms $(S,\NNN)\to (\Spec \OOO,\MMM)$
and whose morphisms are strict $(\Spec \OOO,\MMM)$-morphisms 
between them
(cf. Section 2.3).
By Proposition~\ref{alg},
the stack $([\Spec k[F]/G]\times_{k}\GG_m^l)_{/\scriptsize{\Spec} k[P]\times_{k}\GG_m^l}$ represents
the stack whose objects over $S\to \Spec k[P]\times_{k}\GG_m^l$
are
$\Sigma^0$-FR morphisms  $(S,\NNN)\to (\Spec k[P]\times_{k}\GG_m^l,\NNN_P)$ 
and whose morphisms are strict
$(\Spec k[P]\times_{k}\GG_m^l,\NNN_P)$-morphisms between them.
Here we abuse notation and write $\NNN_P$
for the log structure associated to 
the natural map
$P\to k[P]\otimes_k\Gamma(\GG_m^l,\OOO_{\GG_m^l})$
(i.e., the canonical log structure on $X_{\Sigma}=\Spec k[P]\times_{k}\GG_m^l$),
and $G$ acts on $\Spec k[F]$
in the same way as above.
Consider the cartesian diagram
\[
\xymatrix{
[\Spec k\{F,t_1,\ldots,t_l\}/G]\ar[d] \ar[r] & [\Spec k[F]/G]\times_{k}\GG_m^l \ar[d]\\
\Spec k\{ P,t_1,\ldots,t_l\}  \ar[r]& \Spec k[P]\times_{k}\GG_m^l, \\
}
\]
where the lower horizontal arrow is
$\alpha:\Spec \OOO\to \Spec k[P]\times_{k}\GG_m^l$.
Then this diagram implies our assertion and we conclude that
$\Phi_{\bar{x}}$ is an isomorphism.

\noindent
(Step 4)  Finally, we shall show that the diagram
\[
\xymatrix{
\XX\times \GG_m^d \ar[d]^{\Phi\times\Phi_0}\ar[r]^{m} & \XX \ar[d]^{\Phi}\\
\XXX_{(\Sigma,\Sigma^0)}\times\Spec k[M]\ar[r]^(0.7){a_{(\Sigma,\Sigma^0)}} & \XXX_{(\Sigma,\Sigma^0)} \\
}
\]
commutes. 
Let $\Psi:\XXX_{(\Sigma,\Sigma^0)}\times\Spec k[M]\to\XX\times \GG_m^d$
be a functor such that $(\Phi\times\Phi_0)\circ \Psi\cong\textup{Id}$
and $\Psi\circ (\Phi\times\Phi_0)\cong\textup{Id}$.
Notice that both $\Phi\circ m\circ \Psi$ and $a_{(\Sigma,\Sigma^0)}$ are liftings of the torus action $X_{\Sigma}\times\Spec k[M]\to X_{\Sigma}$.
Then we have $\Phi\circ m\circ \Psi\cong a_{(\Sigma,\Sigma^0)}$ because
a lifting (as a functor)
$\XXX_{(\Sigma,\Sigma^0)}\times\Spec k[M]\to \XXX_{(\Sigma,\Sigma^0)}$ of the torus action $X_{\Sigma}\times\Spec k[M]\to X_{\Sigma}$ is unique up to a unique isomorphism (cf. Corollary~\ref{uniq}). Thus $\Phi\circ m\cong a_{(\Sigma,\Sigma^0)}\circ (\Phi\times\Phi_0)$.
This completes the proof of Theorem~\ref{main2}.
\QED

{\it Proof of Theorem~\ref{main3}.}
It follows from Theorem~\ref{main} and Theorem~\ref{main2}.
\QED

Let $\XXX$ be an algebraic stack. For a point $a:\Spec K\to \XXX$ with an algebraically closed field $K$ the stabilizer group scheme 
is defined to be $\textup{pr}_1:\Spec K\times_{(a,a),\XXX\times\XXX,\Delta}\XXX\to \Spec K$,
where $\Delta$ is diagonal.
If $\XXX$ is Deligne-Mumford, then the stabilizer group scheme is a finite
group.
The proof of Theorem~\ref{main2} immediately implies:

\begin{Corollary}
\label{levelchar}
Let $\XXX$ be a smooth Deligne-Mumford stack separated and of finite type
over $k$.
Suppose that there exists a coarse moduli map
$\pi:\XXX\to X_{\Sigma}$ to a toric variety
such that $\pi$ is an isomorphism over $T_{\Sigma}$.
Let $V(\rho)$ denote the torus-invariant divisor corresponding to each ray
$\rho$,
and suppose that the order of stabilizer group
of the generic point on $\pi^{-1}(V(\rho))$ is $n_{\rho}$.
Then there exists an isomorphism $\XXX\to \XXX_{(\Sigma,\Sigma^0)}$
over $X_{\Sigma}$, where the level of $\Sigma^0$ on $\rho$
is $n_{\rho}$ for each $\rho$.

\end{Corollary}

\Proof
By the proof of Theorem~\ref{main2},
there exist some stacky fan $(\Sigma,\Sigma^0)$
and an isomorphism $\XXX\cong\XXX_{(\Sigma,\Sigma^0)}$ over $X_{\Sigma}$.
Moreover if the level of $\Sigma^0$ on $\rho$ is $n$,
then by \cite[Proposition 4.13]{I2} the stabilizer group
of the generic point on the torus-invariant divisor on $\XXX_{(\Sigma,\Sigma^0)}$ corresponding to $\rho$ is of the form $\mu_n=\Spec K[X]/(X^n-1)$.
Therefore our claim follows.
\QED

\begin{Remark}
\label{semifinalrem}
In virtue of Theorem~\ref{main2},
one can handle toric triples, regardless of their constructions,
by machinery of toric algebraic stacks \cite{I2} and various approaches.
(See Section 5.)

One reasonable generalization of toric triple to positive characteristics
might be 
a smooth tame Artin stack with finite diagonal that is of finite type
over an algebraically closed field, satisfying (i), (ii), (iii)
in Introduction.
(For the definition of tameness, see \cite{AOV}.
Since the stabilizer group of each point on a toric algebraic stack
is diagonalizable, thus every toric algebraic stack is a tame Artin stack.)
Indeed, toric algebraic stacks defined in \cite{I2} are toric triples
in this sense in arbitrary characteristics.
We conjecture that the geometric characterization theorem holds also in positive characteristics.
\end{Remark}

\begin{Remark}
\label{pict}
Let us denote by $\mathfrak{Torst}$ the 2-category of toric algebraic stacks,
or equivalently (by Theorem~\ref{main3}) 2-category of toric triples
(cf. Section 1).
Let us denote by $\textup{Smtoric}$ (resp. $\textup{Simtoric}$)
the category of smooth (non-singular) toric varieties
(resp. simplicial toric varieties),
whose morphisms are torus-equivariant.
From the results that we have obtained so far,
we have the following commutative diagram (picture).
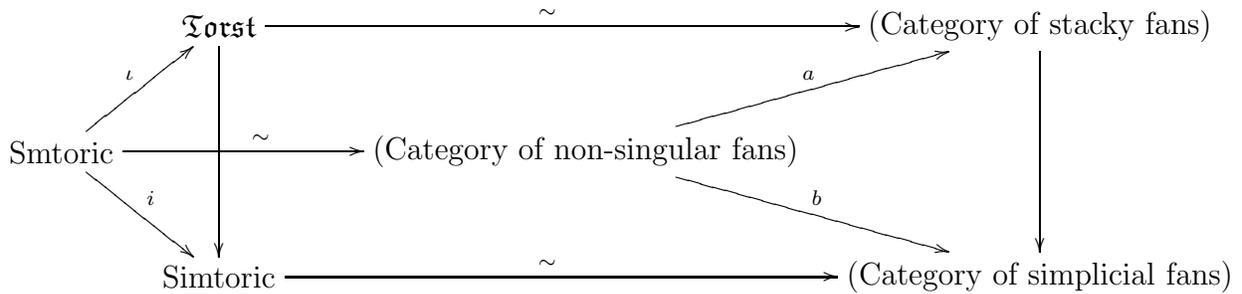
\[
\xymatrix@R=10mm @C=4mm{
 & \mathfrak{Torst}\ar[rrr]^(0.4){\sim} \ar[dd]& & &(\textup{Category\ of\ stacky\ fans}) \ar[dd] \\
\textup{Smtoric} \ar[rd]^{i}\ar[ru]^{\iota}\ar[rrr]^(0.38){\sim}& & & (\textup{Category\ of\ non-singular\ fans})\ar[ru]^{a}\ar[rd]^{b}& \\   &  \textup{Simtoric} \ar[rrr]^(0.4){\sim}& & & (\textup{Category\ of\ simplicial\ fans})\\}
\]
where
$a(\Sigma)=(\Sigma,\Sigma^0_{\extup})$
and $b(\Sigma)=\Sigma$ for a non-singular fan $\Sigma$.
The functors $\iota$ and $i$ are natural inclusion functors.
The functors $\iota$, $i$, $a$ and $b$ are {\it fully faithful}.
All horizontal arrows are {\it equivalences}.
\end{Remark}

\section{Related works}
In this Section we discuss the relationship with \cite{BCS}, \cite{FMN} and \cite{P}.
We work over the complex number field $\CC$.
If no confusion seems to likely arise, we refer to toric triples as
toric stacks.

We first recall the stacky fans introduced in \cite{BCS}.
Let $N$ be a finitely generated abelian group.
Let $\Sigma$ be a simplicial fan in $N\otimes_{\ZZ}\QQ$.
Let $\bar{N}$ be the lattice, that is, the image of $N\to N\otimes_{\ZZ}\QQ$.
For any $b\in N$, we denote by $\bar{b}$ the image of $b$ in $\bar{N}$.
Let $\{\rho_1,\ldots,\rho_r\}$ be the set of rays of $\Sigma$.
Let $\{b_1,\ldots,b_r\}$ be the set of elements of $N$
such that each $\bar{b}_i$ spans $\rho_i$.
The set $\{b_1,\ldots,b_r\}$ gives rise to the
homomorphism $\beta:\ZZ^r\to N$.
The triple $\SIGMA=(N,\Sigma,\beta:\ZZ^r\to N)$ is called a stacky fan.
If $N$ is free, then we say that $\SIGMA$ is reduced.
Every stacky fan $\SIGMA$ has the natural underlying reduced stacky fan
$\SIGMA_{\xtup}=(\bar{N},\Sigma, \bar{\beta}:\ZZ^r\to \bar{N})$,
where $\bar{N}=N/(\textup{torsion})$ is the lattice and
$\bar{\beta}$ is defined to be the composite $\ZZ^r\to N\to \bar{N}$.

Let $\SIGMA=(N,\Sigma,\beta)$ be a reduced stacky fan.
Let $\beta_{\ge0}:\ZZ_{\ge 0}^r\to N$ be the map induced by the restriction
of $\beta$.
The intersection $\beta_{\ge 0}(\ZZ^r_{\ge 0})\cap|\Sigma|$
forms a free-net of $\Sigma$.
The pair $(\Sigma,\beta_{\ge 0}(\ZZ^r_{\ge 0})\cap|\Sigma|)$
is a stacky fan in the sense of Definition~\ref{stfan}.
Conversely, every stacky fan $(\Sigma,\Sigma^0)$ in Definition~\ref{stfan}
is obtained from a unique reduced stacky fan $\SIGMA=(N,\Sigma,\beta)$.
It gives rise to a one-to-one bijective correspondence between
reduced stacky fans and stack fans in the sense of Definition~\ref{stfan}.
In order to avoid confusions, in this Section we refer
to a stacky fan in the sense of Definition~\ref{stfan}
as a framed stacky fan.

Let $\SIGMA=(N,\Sigma,\beta)$ be a stacky fan.
Assume that the rays span the vector space $N\otimes_{\ZZ}\QQ$.
In \cite{BCS}, modelling the construction
of D. Cox \cite{C}, toric Deligne-Mumford stack $\XX(\SIGMA)$ is 
constructed as a quotient stack $[Z/G]$.
There exists a coarse moduli map $\pi(\SIGMA):\XX(\SIGMA)\to X_{\Sigma}$.

Suppose that $\SIGMA$ is reduced and let $(\Sigma,\Sigma^0)$
be the corresponding framed stacky fan.
Then we have:

\begin{Proposition}
\label{vsBCS}
There exists an isomorphism
\[
\XXX_{(\Sigma,\Sigma^0)}\stackrel{\sim}{\longrightarrow}\XX(\SIGMA)
\]
of algebraic stacks over $X_{\Sigma}$.
\end{Proposition}

\Proof
We will prove this Proposition by applying the geometric characterization
Theorem~\ref{main2} and Corollary~\ref{levelchar}.
Let $d$ be the rank of $N$.
The stack $\XX(\SIGMA)$ is a smooth $d$-dimensional Deligne-Mumford stack separated
and of finite type over $\CC$, and
its coarse moduli space is the toric variety
$X_{\Sigma}$ (see \cite[Lemma 3.1, Proposition 3.2 and Proposition 3.7]{BCS}).
From the quotient construction, $\XX(\SIGMA)$
has a torus embedding $\GG_m^d\to \XX(\SIGMA)$.
By Corollary~\ref{levelchar}, to prove our Proposition
it suffices to check that the order of stabilizer group at the generic point
of
$\pi(\SIGMA)^{-1}(V(\rho_i))$ is equal to the level of $\Sigma^0$ on $\rho_i$.
Here $V(\rho_i)$ is the torus-invariant divisor corresponding to
$\rho_i$. To this end, we may assume that $\Sigma$ is a complete fan.
Then \cite[Proposition 4.7]{BCS} implies that
the order of stabilizer group at the generic point of
$\pi(\SIGMA)^{-1}(V(\rho_i))$
is the level $n_{\rho_i}$ of $\rho_i$.
\QED

\begin{Remark}
The explicit construction of $\XX(\SIGMA)$
plays no essential role in the proof of Proposition~\ref{vsBCS},
and the proof uses only some intrinsic properties.
It may show the flexibility of our results.
If a new approach (construction) to this subject is proposed
(in the future),
then the category of toric triples will provide a useful bridge.
Taking it into account, we believe that a good attitude is to have various approaches
at one's disposal and to feel free in choosing one of them
depending on situations.
\end{Remark}

Let $\SIGMA=(N=N'\oplus \ZZ/w_1\ZZ\oplus\cdots \oplus \ZZ/w_t\ZZ,\Sigma,\beta)$ be a stacky fan such that $N'$ is a free abelian group,
and let $\SIGMA_{\xtup}$ be the associated reduced
stacky fan. There is a morphism
$\XX(\SIGMA)\to \XX(\SIGMA_{\xtup})$,
which is a finite abelian gerbe.
This structure is obtained by a simple technique
called ``taking $n$-th roots of an invertible sheaf".
We will explain it.
Recall the notion of the stack of roots of an invertible sheaf
(for example, see \cite{Cad}). 
Let $X$ be an algebraic stack
and $\LLL$ an invertible sheaf on $X$.
Let $P:X\to B\GG_m$ be the morphism to
the classifying stack of $\GG_m$, that corresponds to $\LLL$.
Let $l$ be a positive integer and
$f_l:B\GG_m\to B\GG_m$ the morphism associated to
$l:\GG_m\to \GG_m:\ g\mapsto g^l$.
Then the stack of $l$-th roots of $\LLL$
is defined to be $X\times_{P,B\GG_m,f_l}B\GG_m$.
The reason why this stack is called the stacks of $l$-th roots
is that it has the following modular interpretation:
Objects of $X\times_{P,B\GG_m,f_l}B\GG_m$ over a scheme $S$
are triples $(S\to X, \MMM,\phi: \MMM^{\otimes l}\to \LLL)$,
where $\MMM$ is an invertible sheaf on $S$ and $\phi$ is
an isomorphism. A morphism of triples is defined in a natural manner.
The first projection $X\times_{P,B\GG_m,f_l}B\GG_m\to X$
forgets data $\MMM$ and $\phi$.
We write $X(\LLL^{1/l})$ for the stacks of $l$-th roots of $\LLL$.
In \cite[Proposition 2.9, Remark 2.10]{JT} or \cite[Proposition 3.1]{P}
it was observed and shown that
$\XX(\SIGMA)$ is a finite abelian gerbe
over $\XX(\SIGMA_{\xtup})$ which is obtained
by using the stacks of roots of invertible sheaves.
In the light of Proposition~\ref{vsBCS},
it is stated as follows:

\begin{Corollary}
\label{vsexplicit}
Let $(\Sigma,\Sigma^0)$ be the framed stacky fan
that corresponds to $\SIGMA_{\xtup}$.
Let $b_{i,j}\in \ZZ/w_j\ZZ$ be the image of $b_i \in N$ in $\ZZ/w_j\ZZ$.
$($We may regard $b_{i,j}$ as an element of $\{0,\ldots,w_{j}-1\}$.$)$
Let $\NNN_{i}$ be an invertible sheaf on $\XXX_{(\Sigma,\Sigma^0)}$,
which is associated to the torus-invariant divisor
corresponding to $\rho_i$.
Let $\LLL_{j}=\otimes_{i}\NNN_i^{\otimes b_{i,j}}$.
Then $\XX(\SIGMA)$ is isomorphic to
\[
\XXX_{(\Sigma,\Sigma^0)}(\LLL_1^{1/w_1})\times_{\XXX_{(\Sigma,\Sigma^0)}}\cdots \times_{\XXX_{(\Sigma,\Sigma^0)}}\XXX_{(\Sigma,\Sigma^0)}(\LLL_r^{1/w_r}).
\]

\end{Corollary}

It is known that every separated normal Deligne-Mumford stack
is a gerbe over a Deligne-Mumford stack that is generically a scheme.
We will consider an intrinsic characterization of
toric Deligne-Mumford stacks in the sense of \cite{BCS}
from the viewpoint of gerbes.
Since the construction in \cite{BCS}
employed
the idea of Cox, thus we need to impose the assumption
that the rays $\{\rho_1,\ldots,\rho_r\}$
span the vector space $N\otimes_{\ZZ}\QQ$.
In order to fit in with \cite{BCS},
we consider the following condition of toric stacks (toric triples).
A toric stack (triple) $\XX$ is said to be full
if $\XX$ has no splitting $\XX'\times \GG_m^p$,
such that $\XX'$ is a toric stack and $p$ is a positive integer.

Let $\XXX$ be an algebraic stack.
We say that an algebraic stack $\YYY\to \XXX$ is a {\it polyroots gerbe} over $\XXX$ if it
has the form of the composite
\[
\XXX(\LLL_1^{1/n_1})(\LLL_2^{1/n_2})\cdots (\LLL_k^{1/n_k})\to \XXX(\LLL_1^{1/n_1})(\LLL_2^{1/n_2})\cdots (\LLL_{k-1}^{1/n_{k-1}})\to \cdots\to \XXX,
\]
where 
$\LLL_i$ is an invertible sheaf on
$\XXX(\LLL_1^{1/n_1})(\LLL_2^{1/n_2})\cdots (\LLL_{i-1}^{1/n_{i-1}})$
for $i\ge 2$
and $\LLL_1$ is an invertible sheaf on $\XXX$.
For example, if $\LLL_a$ and $\LLL_b$ is invertible sheaves on $\XXX$,
and $\YYY=\XXX(\LLL_a^{1/n})\times_{\XXX}\XXX(\LLL_b^{1/n'})$,
then $\YYY$ is the composite $\XXX(\LLL_a^{1/n})(\LLL_b|_{\XXX(\LLL_a)}^{1/n'})\to \XXX(\LLL_a)\to \XXX$, and thus it is a polyroots gerbe over $\XXX$.

\begin{Proposition}
A toric Deligne-Mumford stack in the sense of \cite{BCS}
is precisely characterized as
a polyroots gerbe over a full toric stack $($i.e., full toric triple in our sense$)$.
\end{Proposition}

\Proof
Note first that by Corollary~\ref{vsexplicit}
every toric Deligne-Mumford stack $\XX(\SIGMA)$
is a polyroots gerbe over some toric stack $\XXX_{(\Sigma,\Sigma^0)}$,
so $\XX(\SIGMA)$ is a polyroots gerbe over $\XXX_{(\Sigma,\Sigma^0)}$.
Since the condition on $\Sigma$ that rays span the vector space $N\otimes_{\ZZ}\QQ$
is equivalent to the condition that $\XXX_{(\Sigma,\Sigma^0)}$
is full, thus
every toric Deligne-Mumford stack $\XX(\SIGMA)$
is a polyroots gerbe over a full toric stack (triple).
Thus we will prove the converse.
It suffices to prove that for any toric Deligne-Mumford stack $\XX(\SIGMA)$ and any invertible sheaf
$\LLL$ on it, the stack $\XX(\SIGMA)(\LLL^{1/n})$
is also a toric Deligne-Mumford stack in the sense of \cite{BCS}.
Here $n$ is a positive integer.
Note that $\XX(\SIGMA)$ is the quotient stack $[Z/G]$
where $G$ is a diagonalizable group and $Z$ is an open subset
of an affine space $\mathbb{A}^q$ such that codimension of
the complement $\mathbb{A}^q-Z$
is greater than 1.
Thus the Picard group $\mathbb{A}^q$ is naturally isomorphic to
that of $Z$, so every invertible sheaf on $Z$ is trivial,
that is, every principal $\GG_m$-bundle on $Z$ is trivial.
Therefore every principal $\GG_m$-bundle on $[Z/G]$
has the form $[Z\times \GG_m/G]\to[Z/G]$
where the action of $G$ on $Z\times \GG_m$ arises from the action of $G$ on $Z$
and some character $\lambda:G\to \GG_m$.
Let $\alpha:[Z/G]\to B\GG_m$ be the morphism
induced by $Z\to \Spec \CC$ and $\lambda:G\to \GG_m$.
Notice that $\alpha:[Z/G]\to B\GG_m$
is the composite $[Z/G]\to BG\to B\GG_m$, where the first morphism is
induced by the $G$-equivariant morphism $Z\to \Spec \CC$
and the second morphism is induced by $\lambda:G\to \GG_m$.
Let $l:\GG_m\to \GG_m:\ g\mapsto g^l$
and let $\tilde{G}:=G\times_{\lambda,\GG_m,l}\GG_m$.
Then we obtain
the diagram
\[
\xymatrix{
\mathcal{G} \ar[r]\ar[d]& B\tilde{G} \ar[d]\ar[r] &B\GG_m\ar[d]^{f_l}\\
[Z/G] \ar[r]\ & BG  \ar[r] & B\GG_m
 }
\]
where the left square is a cartesian diagram and the right vertical
morphism is induced by $l:\GG_m\to \GG_m$.
Then
by \cite[Corollary 1.2]{J} $\mathcal{G}$ is a toric Deligne-Mumford stack.
Since $[Z\times \GG_m/G]\cong [Z/G]\times_{\alpha,B\GG_m}\Spec \CC$,
the morphism
$\alpha$ corresponds to the principal $\GG_m$-bundle
$[Z\times\GG_m/G]\to [Z/G]$.
Thus if $\mathcal{G}\cong [Z/G]\times_{\alpha,B\GG_m,f_l}B\GG_m$,
it follows our claim. Thus it suffices to check
that the right square is a cartesian diagram.
Indeed, there exists a natural isomorphism
$B\mu_l \cong \Spec \CC\times_{B\GG_m,l}B\GG_m$
and thus we have $B\mu_l \cong \Spec \CC\times_{BG}(BG\times_{B\GG_m,f_l}B\GG_m)$.
Moreover there exists a natural isomorphism
$B\mu_l \cong \Spec \CC\times_{BG}B\tilde{G}$
because the kernel of $\tilde{G}\to G$ is $\mu_l$.
Consider the natural morphism
$B\tilde{G}\to BG\times_{B\GG_m,f_l}B\GG_m$ over $BG$.
Its pullback by the flat surjective morphism $\Spec \CC \to BG$
is an isomorphism $B\mu_l\to B\mu_l$.
Thus $B\tilde{G}\cong  BG\times_{B\GG_m,f_l}B\GG_m$.
Hence our proof completes.
\QED

{\it Relation to} \cite{FMN}. In \cite{FMN}, Fantechi-Mann-Nironi generalize
the notion of toric triples introduced in this paper.
In order to fit in with
the framework of \cite{BCS},
they introduced ``DM torus'' which is a torus with a trivial gerbe
structure and considered
actions of DM tori on algebraic stacks.
Following the point of view that ``toric objects" should be characterized
by torus embeddings and actions,
they discuss a geometric
characterization of toric Deligne-Mumford stack in th sense of \cite{BCS}
by means of smooth Deligne-Mumford stacks with DM torus embeddings
and actions (cf. \cite[Theorem II]{FMN}). Namely, the embeddings and actions
of DM tori
provide gerbe structures
on toric triples, discussed above.

\vspace{1mm}

{\it Notes.}
The former version of this paper was posted on arXiv server during
 December 2006
and in it the main results of this paper were proven,
whereas \cite{FMN} appeared on arXiv in August 2007.

\vspace{2mm}

{\it Relation to} \cite{P}.
In \cite{P}, Perroni studied 2-isomorphism classes
of all 1-morphisms between toric Deligne-Mumford stacks
in the sense of \cite{BCS}. The method and description
are parallel to \cite[section 3]{C2}.
Let $\XX(\SIGMA)$ and $\XX(\SIGMA')$ be toric Deligne-Mumford stacks.
Suppose that $\XX(\SIGMA)$ is proper over $\CC$.
Then Perroni gave a description of 2-isomorphism classes of 1-morphisms from
$\XX(\SIGMA)$ to $\XX(\SIGMA')$
in terms of homogeneous polynomials of $\XX(\SIGMA)$.
(For details, see \cite[Section 5]{P}.)
Assume that $\SIGMA$ and $\SIGMA'$ are reduced.
If the morphism $\tilde{f}:\XX(\SIGMA)\to \XX(\SIGMA')$
associated to a system of homogeneous polynomials
(cf. \cite[Theorem 5.1]{C2}) induces a
torus-equivariant morphism $f:X_{\Sigma}\to X_{\Sigma'}$,
then by Corollary~\ref{uniq} and~\ref{equivcha},
$\tilde{f}$ is a torus-equivariant morphism.
Namely, if $(\Sigma,\Sigma^0)$ and $(\Sigma',\Sigma^{'0})$
denote framed stacky fans corresponding to
$\SIGMA$ and $\SIGMA'$ respectively,
then the morphism $\Sigma\to \Sigma'$ of fans corresponding to $f:X_{\Sigma}\to X_{\Sigma'}$
induces $(\Sigma,\Sigma^0)\to (\Sigma',\Sigma^{'0})$,
and through isomorphisms $\XXX_{(\Sigma,\Sigma^0)}\cong \XX(\SIGMA)$ and $\XXX_{(\Sigma',\Sigma^{'0})}\cong\XX(\SIGMA')$, the morphism $\XXX_{(\Sigma,\Sigma^0)}\to \XXX_{(\Sigma',\Sigma^{'0})}$ associated to $(\Sigma,\Sigma^0)\to (\Sigma',\Sigma^{'0})$ is identified with $\tilde{f}$.

\end{document}